\documentclass[11pt,a4paper,english,reqno]{amsart}
\usepackage{amsmath,amssymb,amsfonts,epsfig,mathrsfs}
\usepackage[T1]{fontenc}

\makeatletter
\@namedef{subjclassname@2020}{%
  \textup{2020} Mathematics Subject Classification}
\makeatother

\usepackage{color}
\usepackage{array}
\usepackage{amsthm}
\usepackage{amstext}
\usepackage{amsmath,stackengine}
\usepackage{graphicx}
\usepackage{setspace}
\usepackage{comment}
\usepackage[margin=2.5cm]{geometry}
\usepackage{bbm}
\usepackage{color}
\usepackage{enumitem}
\usepackage{undertilde}
\allowdisplaybreaks[4]

\usepackage{pgfplots}

\usepackage{appendix}

\usepackage{amscd,psfrag}
\usepackage{yhmath}
\usepackage[mathscr]{eucal}

\usepackage{slashed}

\makeatletter
\pdfpageheight\paperheight
\pdfpagewidth\paperwidth

\newcommand{\bra}{\left\langle}
\newcommand{\ket}{\right\rangle}

\usepackage{mathrsfs}

\usepackage{epstopdf}
\usepackage{chngcntr}
\counterwithin{figure}{section}
\usepackage{mathrsfs}

\usepackage{indentfirst}	

\usepackage[normalem]{ulem}
\theoremstyle{plain}

\newtheorem{definition}{Definition}[section]
\newtheorem{theorem}[definition]{Theorem}
\newtheorem*{theorem*}{Theorem}

\newtheorem{remark}[definition]{Remark}

\newtheorem*{remark*}{Remark}
\newtheorem*{sideremark*}{Side Remark}\newtheorem*{mt*}{Main Theorem}

\newtheorem*{claim*}{Claim}
\newtheorem*{q*}{Question}
\newtheorem{lemma}[definition]{Lemma}

\newtheorem*{corollary*}{Corollary}
\newtheorem*{proposition*}{Proposition}
\newtheorem{example}[definition]{Example}
\newtheorem{proposition}[definition]{Proposition}

\newcommand{\R}{\mathbb{R}}

\newcommand{\na}{\nabla}

\newcommand{\dd}{{\rm d}}
\newcommand{\p}{\partial}
\newcommand{\e}{\epsilon}
\newcommand{\emb}{\hookrightarrow}
\newcommand{\weak}{\rightharpoonup}
\newcommand{\map}{\rightarrow}
\newcommand{\G}{\Gamma}

\newcommand{\dvg}{\dd V_g}

\newcommand{\lh}{\overline{h}}

\newcommand{\mean}{\mathfrak{m}}
\newcommand{\ric}{{\rm Ric}}
\newcommand{\scal}{{\rm scal}}
\newcommand{\gq}{{g^{(q)}}}
\newcommand{\kq}{{\kappa^{(q)}}}\newcommand{\kqone}{{\kappa^{(q)}_1}}\newcommand{\kqtwo}{{\kappa^{(q)}_2}}\newcommand{\hq}{{h^{(q)}}}\newcommand{\Hq}{{H^{(q)}}}
\newcommand{\gaq}{{\gamma^{(q)}}}

\newcommand{\wstar}{\mathrel{\ensurestackMath{\stackon[1pt]{\rightharpoonup}{\scriptstyle\ast}}}}

\newcommand{\ttwo}{{\mathbb{T}^2}}

\newcommand{\tor}{{\mathbf{T}(a,b)}}

\newcommand{\torc}{{\mathbf{T}_c(a,b)}}

\allowdisplaybreaks[4]

\def\XXint#1#2#3{{\setbox0=\hbox{$#1{#2#3}{\int}$ }
\vcenter{\hbox{$#2#3$ }}\kern-.6\wd0}}

\numberwithin{equation}{section}
\numberwithin{figure}{section}

\title{From the Nash--Kuiper Theorem of Isometric Embeddings to the Euler Equations for Steady Fluid Motions: Analogues, Examples, and Extensions}

\author{Siran Li}
\address{Siran Li: School of Mathematical Sciences $\&$ IMA-Shanghai, Shanghai Jiao Tong University, No.~6 Science Buildings,
800 Dongchuan Road, Minhang District, Shanghai, China (200240)}

\email{\texttt{siran.li@sjtu.edu.cn}}

\author{Marshall Slemrod}
\address{Marshall Slemrod: Department of Mathematics, University of Wisconsin, Madison, Madison, WI 53706, USA.}
\email{\texttt{slemrod@math.wisc.edu}}


\date{\today}

\begin{document}

\maketitle

\begin{center}
\emph{This paper is dedicated to our friend and mentor\\ Costas Dafermos on the occasion of his $81^{\text{st}}$ birthday}
\end{center}

\begin{abstract}
Direct linkages between regular or irregular isometric embeddings of surfaces and  steady compressible or incompressible fluid dynamics are investigated in this paper. For a surface $(M,g)$ isometrically embedded in $\R^3$, we construct a mapping which sends the second fundamental form of the embedding to the density, velocity, and pressure of steady fluid flows on $(M,g)$. From the PDE perspectives, this mapping sends solutions to the Gauss--Codazzi equations to the steady Euler equations. Several families of special solutions of physical or geometrical significance are studied in detail, including the Chaplygin gas on standard and flat tori, as well as the irregular isometric embeddings of the flat torus. We also discuss tentative  extensions to multi-dimensions.
\end{abstract}

\section{Introduction}
In recent years, an important research project attracting enormous attention in PDE, differential geometry, and mathematical hydrodynamics communities is concerned with irregular isometric embeddings of Riemannian manifolds and wild solutions for fluid dynamical PDEs. A sequence of papers  has appeared as an outgrowth of the celebrated results of Nash and Kuiper \cite{nash, k1, k2} on $C^1$-isometric embeddings of  Riemannian manifolds into Euclidean spaces. These works are devoted both to sharpening the regularity of embeddings in the Nash-Kuiper theorem (\cite{cds, x}), and to proving non-uniqueness of solutions for problems arising in continuum mechanics, especially the Euler equations for compressible and incompressible fluid flows (\cite{a1, a2, a3}, amongst many others). See De Lellis--Sz\'{e}kelyhidi Jr. \cite{bams} for a thorough survey.

In these works, in-depth, far-reaching connections between geometry and fluid --- more precisely, between the isometric embedding problem and the fluid dynamical PDE, \emph{e.g.}, the Euler equations --- have been established and explored. Such connections are better  understood on the \emph{methodological} level. Indeed, a central technique in common for the above works is convex integration, which enables researchers to  establish, in a constructive manner, the existence of irregular isometric embeddings or the non-uniqueness of weak solutions to fluid dynamical PDE. Systematization and further developments of convex integration     culminated in Gromov's $h$-principle and partial differential relations, which have become important theories in geometry and topology today. See \cite{gr1, gr2}.

The above theories and techniques have also led to important advancements in mathematical hydrodynamics and analysis of PDEs in continuum mechanics. One cornerstone is the final resolution of the Onsager conjecture for incompressible Euler equations. See \cite{a1, a2, bams, a3} by Buckmaster, De Lellis, Isett, and Sz\'{e}kelyhidi Jr., and the many references cited therein. Let us also mention that recently, complementing the above theoretical developments are the breakthrough
computational results of Borrelli--Jabrane--Lazarus--Thibert \cite{b1, b2, b3}, in which for the first time the convex integration procedure was numerically implemented to produce elegant illustrations
of irregular isometric embeddings.

Despite the aforementioned rapid developments in this field of research, the connections between isometric embeddings and fluid dynamics still await better understanding, beyond the level of  structural similarities in relevant PDE or shared analytical tools (\emph{i.e.}, convex integration). Indeed, to the best of the authors' knowledge, in each instance of its applications 
to continuum mechanics in the literature, the convex integration procedure is carried out on a case by case basis (similarly for the numerical results in \cite{b2, b3}); moreover, the geometrical characteristics of convex integration are, to some extent,  neglected or suppressed. One exception is the recent work \cite{th} by Theilli\`{e}re, in which a unified, systematic convex integration framework has been proposed, which encompasses several important geometric examples including the $C^1$-Nash embeddings.

In this work, we present some investigations, though of a tentative and preliminary nature, on the ``geometry-fluid correspondence'' between the isometric embedding problem and fluid dynamical PDE. Our aim is to establish possible \emph{direct fluid dynamic analogues of isometric embeddings}. That is, we discuss ways of directly translating the isometric embeddings,  regular or irregular, into suitable PDE for fluid motions.

Our tentative, partial answers to this problem can be summarised in the sequel:
\begin{enumerate}
\item
For the convex integration
procedure utilised to isometrically embed a surface $(M,g)$ into $\R^3$, each member of the approximating sequence can be identified with a
solution of the Euler equations for  compressible steady fluid flow on $(M,g)$.
\item
Then, upon {\em ``renormalisation''}, the Nash--Kuiper limit of approximate solutions can be translated into solutions to the steady compressible Euler equation on $(M,g)$.
\item
Of particular interest is the sequence of rescaled standard tori $\torc$ (see \eqref{torc}), whose metrics converge in $L^\infty$ to the flat metric, in a manner reminiscent of the ``Nash wrinkles'' in \cite{nash}. On each member of the approximating sequence of manifolds, as well as on the limiting manifold, \emph{i.e.}, the flat torus, there are special solutions which have both natural geometrical meaning (geodesics on tori) and physical (Chaplygin gas)  significance. 
\end{enumerate}


\smallskip

The remaining parts of the paper are organised as follows.

 First, \S\ref{sec: prelim} contains  background materials on the isometric embedding problem. \S\ref{sec: 2d} presents a link between smooth  isometrically embeddings of a surface $(M,g)$ into $\R^3$ and solutions to the steady compressible Euler equations on $(M,g)$. In \S\ref{sec: examples}, two families of special examples for the general theory established in \S\ref{sec: 2d} are presented. Next, in \S\ref{sec: Irrotational Chaplygin gas on the flat torus} we discuss the same fluid dynamical PDE as in earlier sections on the flat torus or on its irregularly isometrically embedded image in $\R^3$. In \S\ref{sec: renorm} we discuss one possible  correspondence between irregular isometric embeddings and steady Euler equations. In particular, we introduce a renormalisation process of successive rescaling and identify the ``renormalised limit'' of Nash--Kuiper iterations as weak solutions to the Euler equations. Finally, some preliminary results on  extensions to manifolds of higher dimensions are reported in \S\ref{sec: multi D}.

Two Appendices~\ref{sec: app A} and \ref{sec: appendix, variation} discuss, respectively, the geometry of rescaled standard tori and study of Chaplygin gas on higher dimensional manifolds. The appendices supplement \S\S\ref{sec: examples} $\&$ \ref{sec: Irrotational Chaplygin gas on the flat torus}.

\section{Geometric Preliminaries}\label{sec: prelim}
This section collects some basics on  Riemannian geometry and the isometric embedding problem. For the moment we restrict our discussions to $2$-dimensional Riemannian manifolds, {\it i.e.}, surfaces.

Let $(M, g)$ be a $2$-dimensional Riemannian manifold, and consider an arbitrary local co-ordinate system $\left\{x^i: i=1,2\right\}$. The distance on $M$ is given by the metric $g$ (\emph{a.k.a.} first fundamental form) via $\dd s^2 = g_{ij} \dd x^i \otimes \dd x^j$. A map $${\bf y}: (M,g) \to \R^3$$ is an isometric embedding if ${\bf y}$ and $\dd {\bf y}$ are both injective and ${\bf y}$ is an isometry:
\begin{equation*}
\sum_{k=1}^3 \frac{\p {\bf y}^k}{\p x^i} \frac{\p {\bf y}^k}{\p x^j} = g_{ij}
\end{equation*}
for each $i,j\in\{1,2\}$. That is, the intrinsic distance on $M$ given by $g$ equals to the ambient Euclidean distance on the image ${\bf y}(M) \subset \R^3$.

Isometric embeddings has been a central topic in the development of differential geometry and nonlinear analysis $\&$ PDE. See Han--Hong \cite{hh} for a thorough introduction.  One crucial problem is concerned with the existence of an  isometric embedding: given a Riemannian manifold $(M,g)$, find an isometric embedding ${\bf y}$ as above. It amounts to solving for the extrinsic geometry  --- in contrast to the intrinsic geometry, namely the geometric quantities determined by $g$. In the case that $M$ is 2-dimensional, the extrinsic geometry  of ${\bf y}: M \to \R^3$ is completely characterised by the {\em second fundamental form} $H = \left\{H_{ij}:1\leq i,j\leq 2\right\}$: 
\begin{equation}
H_{ij} := \frac{\p^2 {\bf y}}{\p x^i \p x^j} \cdot \nu,
\end{equation}
where $\nu$ is the outward unit normal vector field to ${\bf y}(M)$.

There is a well-known necessary condition for the existence of (smooth) isometric embeddings of a surface $(M,g)$ into $\R^3$: the second fundamental form must satisfy the {\em Gauss} and {\em Codazzi} equations: for $i,j,k,l \in \{1,2\}$, it holds that
\begin{eqnarray}
&&R_{ijkl} = H_{ik}H_{jl} - H_{il}H_{jk},\\
&&\na_i H_{jk} = \na_j H_{ik}.\label{xxx}
\end{eqnarray}
Here $\na$ is the covariant derivative associated to the Levi-Civita connection on $(M,g)$. The connection is fully characterised by the Christoffel symbols $$\G^i_{jk} := \frac{1}{2}g^{il}\{\p_j g_{kl} + \p_k g_{jl} - \p_l g_{jk}\},$$ where $g^{-1}=\{g^{ij}\}$. Here and hereafter, we adopt Einstein's summation convention:  repeated lower and upper indices are understood as being summed over. Then, the  Riemann curvature tensor is given by $$R_{lijk} := g_{lp} \left\{\p_j \G^p_{ik} -\p_k \G^p_{ij} + \G^p_{jq} \G^q_{ik} - \G^p_{kq} \G^q_{ij}\right\}.$$ In dimension $2$, the only nontrivial component of the Riemann curvature is $R_{1212}$. One defines the {\em Gauss curvature} by $$\kappa := \frac{R_{1212}}{ \det\,g},$$ which is the only intrinsic curvature for a surface. Let us also mention that $\na g=0$, commonly known as the Ricci identity.

When $M$ is simply-connected, the Gauss--Codazzi equations are also sufficient for the existence of isometric embeddings for $\dim M =2$. This is known as the {\em fundamental theorem of surface theory}; its proof in the case of lower regularity is given by S. Mardare (\cite{ma1, ma2, ma3}).

In the introduction we have discussed the Nash--Kuiper theorem (\cite{nash, k1, k2}). The statement is as follows. Note that an improved version of this theorem has been established by Conti--De Lellis--Sz\'{e}kelyhidi Jr. 
 \cite[Theorem~6.1]{cds}. We shall refer to it in \S\S\ref{sec: Irrotational Chaplygin gas on the flat torus} $\&$ \ref{sec: renorm} below. 
\begin{theorem}\label{thm: nash kuiper}
Let $(M,g)$ be a Riemannian manifold of dimension $n$. Let ${\bf y}_\star : M \map \R^{n+k}$ be a smooth embedding such that $\p_i {\bf y}_\star \cdot \p_j {\bf y}_\star <g_{ij}$ (as quadratic forms) and $k\geq 1$. Then for any $\e>0$ there exists a $C^1$-isometric embedding ${\bf y}: (M,g) \map \R^{n+k}$ such that $\|{\bf y}-{\bf y}_\star\|_{C^0(M)} \leq \e$.
\end{theorem}

Typical applications of the Nash--Kuiper Theorem \ref{thm: nash kuiper} include the construction of a $C^1$-isometric embedding (or even a  $C^{1,\frac{1}{7}-\e}$-isometric embedding) of the flat 2-torus into $\R^3$. Another example is a ``corrugated 2-sphere''; that is, a round 2-sphere of radius $r$ isometrically embedded into $\R^3$ while being $C^0$-close to a sphere of a smller radius $r_0$. A visual representation is given by Bartzos--Borrelli--Denis--Lazarus--Rohmer--Thibert \cite{b1}.

\section{From embedded surfaces to Euler: 2D Smooth Solutions} \label{sec: 2d}

\subsection{Gauss-Codazzi and steady Euler equations}

This section is mainly concerned with the following question:

\begin{quote}
Given a smooth surface $(M,g)$ isometrically embedded in $\R^3$ with second fundamental form $H = \{H_{ij}\}_{1 \leq i,j \leq 2}$, can we identify natural variables $\rho, v, p$ in terms of $g,H$ which represent, respectively, the density, velocity, and pressure of a steady fluid \emph{on $(M,g)$}?
\end{quote}

More precisely, consider the following two PDE systems, one geometric and one physical: 

{\bf 1.} Let $(M,g)$ be a smooth surface isometrically embedded into $\R^3$. Its second fundamental form $H = \{H_{ij}\}_{1 \leq i,j \leq 2}$  solves the Gauss--Codazzi equations: 
\begin{equation}\label{gauss}
H_{ij}H_{kl} - H_{ik} H_{jl} = R_{iljk}
\end{equation}
and
\begin{equation}\label{codazzi}
\na_i H^k_j = \na_j H^k_i.
\end{equation}

{\bf 2.} The Euler conservation laws of mass and momentum for {\em steady} fluid on $(M,g)$:
\begin{equation}\label{E1}
\na_k (\rho v^k) = 0
\end{equation}
and
\begin{equation}\label{E2}
\na_k P^{ik}=0,
\end{equation}
where $P^{ik}$ is the {\em stress-energy tensor}:
\begin{equation}\label{PP}
P^{ik} = \rho v^iv^k + g^{ik}p.
\end{equation}
These are just the usual fluid equations with covariant differentiation; see for example Anco {\it et al} \cite{anco} and Arnold--Khesin \cite{ak}. The idea should not be surprising, since it is building block of the ``matter'' description in general relativity; see Carroll \cite{carroll}.

Our question, formulated in terms of PDE, is as follows: \emph{given a smooth solution $H_{ij}$ to the Gauss--Codazzi system \eqref{gauss} and \eqref{codazzi}, look for fluid variables $\{\rho, v, p\}$ which satisfy the steady Euler equations \eqref{E1} and \eqref{E2} on $(M,g)$}. That is, we hope to provide a direct link from isometric embedding of surfaces to fluid dynamical PDEs. In this section we restrict ourselves to the $C^\infty$-setting; isometric embeddings of weak regularity as in Nash's setting \cite{nash} will be discussed in the subsequent section.

\subsection{Step 1: Identifying the stress-energy tensor}


To begin with, we rewrite via Ricci's identity the balance law of linear momentum \eqref{E2} as\begin{equation}\label{new, E2}
\na_k P^k_j = 0.
\end{equation}
Given the solution $\left\{H^i_j\right\}_{1\leq i,j \leq 2}$ to the Codazzi equation~\eqref{codazzi}, the following choice of $\left\{P^i_j\right\}_{1\leq i,j \leq 2}$ clearly solves the balance of momentum \eqref{new, E2}:
\begin{equation}\label{P def, 2d}
P^1_1 = H^2_2,\quad P^1_2 = -H^1_2, \quad P^2_1=-H^2_1,\quad P^2_2 =H^1_1.
\end{equation}

It is important to note that $\left\{P^i_j\right\}$ can be defined globally on $(M,g)$ as a $(1,1)$-tensor. Indeed, the second fundamental form $H=\{H_{ij}\} \in \G(T^*M \otimes T^*M)$ is a contravariant $2$-tensor, and let us write $H' := g^{-1} H = \left\{H^i_j\right\} \in \G(T^*M \otimes TM)$, a $(1,1)$-tensor obtained by raising one index via contraction with metric $g$. Then 
\begin{align*}
P'=\left\{P^i_j\right\} := {\rm Adj}\, H',
\end{align*}
where the adjugate matrix can be defined globally on $(M,g)$ via multilinear algebra. Note that 
\begin{equation}\label{new, kappa}
\det H' = \det P' =\kappa = \text{the Gauss curvature of $(M,g)$},
\end{equation}
and similarly,
\begin{equation}\label{new, mean}
{\rm tr} \,H' = {\rm tr} \,P' =\mean = \text{the mean curvature of $(M,g) \emb \R^3$}.
\end{equation}
Here and hereafter, $tr$ without subscripts means that the trace is taken with respect to the canonical Euclidean  coordinates.

\subsection{Step 2: Identifying the pressure}

In this step, we shall identify the pressure $p$ with the smaller principal curvature of $(M,g) \emb \R^3$. Let us first introduce some notations: we write $\mean$ for the mean curvature of $M$,  namely that $$\mean = {\rm tr}_g H = H^1_1 + H^2_2.$$ We use the unconventional notation $\mean$ for the mean curvature since $H$ has already been reserved for the second fundamental form. Also note that $\mean := {\rm tr}_g H$ here instead of $\frac{1}{2}{\rm tr}_g H $. In addition, we let $\{\kappa_+, \kappa_-\}$ be the principal curvatures, where $\kappa_+ \geq \kappa_-$ on $M$, and $|\bullet|_g$ denotes the norm of a vector field or $1$-form with respect to the metric; \emph{i.e.}, $$|v|_g^2 := g^{ij}v_iv_j = g_{ij} v^iv^j  \equiv v^iv_i.$$

In the previous step we have defined $P$ as a function of $H$. We now require it to satisfy \eqref{PP}, which is equivalent to $
P_{ij} = \rho v_iv_j + pg_{ij}$. That is, $P = \rho v\otimes v + pg$ as an identity of  contravariant 2-forms. We deduce that
\begin{align}\label{kappa in 2d}
\kappa &= \det (g^{-1}P) \nonumber\\
&= \det\Big[ \rho g^{-1}\cdot (v\otimes v) + p I_2 \Big]\nonumber\\
&= p^2 + {\rm tr}\left[\rho g^{-1}\cdot (v\otimes v)\right] + \det\left[\rho g^{-1}\cdot (v\otimes v)\right]\nonumber\\
&= p\rho |v|^2_g + p^2.
\end{align}
The first line follows from   
\eqref{new, kappa}, the second line holds as $P$ is the stress-energy tensor for steady fluid, the third line is a linear algebra identity, and the last line follows from the definition of norm $|\bullet|_g$ and that $v\otimes v$ is rank-1, hence having zero determinant. Note that \eqref{kappa in 2d} generalises the relations previously obtained for the Euclidean case in Chen--Slemrod--Wang \cite[Equation~(3.2)]{csw}; \emph{cf.}
 also Acharya--Chen--Li--Slemrod--Wang \cite[Eq.~(8.4)]{five}.

On the other hand, we can compute the mean curvature $\mean = {\rm tr}_g\, H = {\rm tr}\, H'$ from \eqref{new, mean}:
\begin{align}\label{new, expression for m}
\mean = H^1_1+H^2_2 =  P^1_1+P^2_2 = \rho |v|^2_g +2p.
\end{align}
Thus, we find from \eqref{new, expression for m} and \eqref{kappa in 2d}  that $\mean$ and $\kappa$ are related by
\begin{equation}\label{eq for mean and gauss curvature, 2d}
p^2 - \mean p + \kappa = 0.
\end{equation}

Let us now consider the principal curvatures $\{\kappa_+, \kappa_-\}$. Then \eqref{eq for mean and gauss curvature, 2d} can be written as $p^2 - (\kappa_+ + \kappa_-)p + \kappa_+\kappa_- = 0$, so the roots are $p=\kappa_+$ or $p=\kappa_-$. We shall select
\begin{equation}\label{p=kappa 2, 2d}
p=\kappa_- = \text{the smaller principal curvature}.
\end{equation}
Indeed, such choice is  compatible with the meaning of $P$ in fluid dynamics; in particular, with the non-negativity of density function $\rho$. To see this, consider the following three cases:
\begin{itemize}
\item
\underline{Case 1: $\kappa_+ >\kappa_-$ and $\kappa_+$ is nonzero}. If in this case $p=\kappa_+$,  then \eqref{kappa in 2d} becomes $p[\rho|v|_g^2 + (p-\kappa_-)]=0$, which is impossible for $\rho \geq 0$. 
\item
\underline{Case 2: $\kappa_+=\kappa_-$}.  Then the choice for $p$ is non-ambiguous.
\item
\underline{Case 3: $\kappa_+ >\kappa_-$ and $\kappa_+$ is vanishing somewhere}. By choosing $p=\kappa_-$, \eqref{kappa in 2d} becomes $p(\rho|v|^2_g +p)=0$, which is permissible as long as $\rho >0$.
\end{itemize}

\subsection{Step 3: Reformulation of balance laws}
With the identification of $P$ and $p$ as in \eqref{P def, 2d} and \eqref{p=kappa 2, 2d}, we shall reformulate in this step the conservation of momentum as a first-order PDE for the $\rho$-variable. In addition, we obtain the PDE for vorticity, which is transported along Lagrangian trajectories on the surface $(M,g)$. This follows from natural considerations in the study of fluid dynamical PDEs and, as will be clear from the later steps, plays a  crucial role in our search for special solutions to the steady Euler equation with geometrical significance.

Following conventions in fluid dynamics, we define the vorticity  $\omega$ on $(M,g)$:
\begin{equation*}
\omega = \na_1 v_2 - \na_2 v_1,
\end{equation*}
which is viewed either as a 2-form or a scalar field, identified via Hodge duality. We also write 
\begin{align*}
q:=|v|_g,\qquad c := \sqrt{p'(\rho)},\qquad M:=q/c
\end{align*}
for the \emph{flow speed}, \emph{sonic speed}, and \emph{Mach number}, respectively.

Recall from \eqref{kappa in 2d}, \eqref{P def, 2d}, and \eqref{p=kappa 2, 2d} that $\kappa = \kappa_+ \kappa_- = p\rho q^2 + p^2$ and $p = \kappa_-$. Thus
\begin{equation}\label{new, flow speed q}
q^2 = \frac{\kappa_+ - \kappa_-}{\rho}.
\end{equation}
Assume that the conservation of mass \eqref{E1} is satisfied. This will reduce the conservation of momentum \eqref{E2} to
\begin{align}\label{mom}
\frac{\rho}{2}v^k\na_k \left(\frac{\kappa_+ - \kappa_-}{\rho}\right) + v^k \na_k \kappa_- =0.
\end{align}
This equation can be written globally on $(M,g)$: 
\begin{equation}\label{mom, global}
\frac{\rho}{2}\na_v \left(\frac{\kappa_+ - \kappa_-}{\rho}\right) + \na_v \kappa_- =0,
\end{equation}
where $\na_v$ is the covariant derivative with respect to $v \in \G(TM)$.

On the other hand, let us consider the vorticity equation. Assuming the conservation of mass \eqref{E1}, we rewrite  the conservation of momentum \eqref{E2} as
\begin{align*}
\rho v^k\na_k v^i + g^{ik}\na_k \kappa_- = 0.
\end{align*}
Lowering the index $i$ (by multiplying  with $g_{i\ell}$ and evaluating separately $\ell=1,2$), one obtains
\begin{equation*}
\begin{cases}
\rho v^1\na_1 v_1 + \rho v^2\na_2 v_1 + \na_1 \kappa_- = 0,\\
\rho v^1\na_1 v^2 + \rho v^2\na_2v_2 + \na_2\kappa_2- = 0,
\end{cases}
\end{equation*}
which is equivalent to
\begin{equation*}
\begin{cases}
\frac{\rho}{2} \na_1 (q^2) - \rho v^2 \omega + \na_1\kappa_- = 0,\\
\rho v^1\omega + \frac{\rho}{2} \na_2(q^2) + \na_2 \kappa_- = 0,
\end{cases}
\end{equation*}
as $q^2 := |v|^2_g = v^1v_1 + v^2v_2$. Multiplying respectively the first and the second equations by $-v^2$ and $v^1$ and adding them up, we have 
\begin{equation}\label{vorticity eq}
\rho q^2 \omega + \frac{\rho}{2} \na_{v^\perp} (q^2) + \na_{v^\perp} \kappa_- = 0,
\end{equation}
where the differential operator $\na_{v^\perp}$ is given in local co-ordinates by 
\begin{equation*}
\na_{v^\perp} = v_1 \na_2 - v_2 \na_1. 
\end{equation*}
This operator is  globally defined on $(M,g)$. Indeed, $$v^\perp = Jv$$ where the $(1,1)$-tensor field $J \in \G(T^*M \otimes TM)$ is the almost complex structure on  $M$.

We observe a direct consequence of the derivations above. The terms involving derivatives in \eqref{mom, global} and \eqref{vorticity eq} are of the same form, in view of the Bernoulli law \eqref{new, flow speed q}. We thus have
\begin{lemma}\label{lem: momentum + irrotational}
With the choice of variables $p=\kappa_-$, $p\rho q^2+p^2 = \kappa$, and $q := |v|_g$, conservation of  momentum and irrotationality are satisfied whenever the density $\rho$ solves the following PDE: 
\begin{align*}
\frac{\rho}{2}\mathcal{L}\left(\frac{\kappa_+ - \kappa_-}{\rho}\right) + \mathcal{L} \kappa_- =0\qquad \text{ with } \mathcal{L} \in \left\{\na_v, \na_{v^\perp}\right\}.
\end{align*}
\end{lemma}

In \S\ref{sec: examples} below, we shall investigate two families of examples by looking for $\rho$ which satisfies the identity in Lemma~\ref{lem: momentum + irrotational} for \emph{every} first-order differential operator $\mathcal{L}$. However, the general procedure in Step~4 below to solve for $v$ and $\rho$ will not make use of this lemma.

\subsection{Step 4: Solving for $v$ and $\rho$}
Now let us discuss how to solve for $v$ and $\rho$ in terms of $g$ and $H$, with the choice of $p$ and $P$ as in earlier Steps~1 $\&$ 2. Our goal is to find a correspondence $$H \rightsquigarrow \left(\rho, v^1, v^2, p\right),$$ such that $H$ satisfies the Gauss--Codazzi equations~\eqref{gauss}, \eqref{codazzi} on $(M,g)$, and $(\rho, v, p)$ satisfies the  conservation of mass~\eqref{E1} and momentum~\eqref{E2}. For this purpose, we have:
\begin{enumerate}
\item
identified $P^{ij} = \rho v^i v^j + g^{ij}p$ as a whole via $\left\{P^i_j\right\} = {\rm Adj}\left\{H^i_j\right\}$. Then the conservation of momentum is automatically verified, thanks to the Codazzi equation~\eqref{codazzi}; and 
\item
expressed the Gauss equation~\eqref{gauss} in the form of \eqref{kappa in 2d}: $\kappa = p\rho q^2 + p^2$. 
\end{enumerate}
In this way, we have obtained four equations for four unknowns:
\begin{equation}\label{four}
\begin{cases}
\rho\left(v^1\right)^2 + g^{11} p = P^{11},\\\rho v^1 v^2 + g^{12} p = P^{12},\\
\rho\left(v^2\right)^2 + g^{22} p = P^{22},\\
p\rho q^2 + p^2 = \kappa.
\end{cases}
\end{equation}
The right-hand sides depend only on the geometric data, which are deemed as given. Moreover, we have also
\begin{itemize}
\item[(3)]
taken the trace of the first three equations in \eqref{four} to obtain \eqref{new, expression for m}; that is, $\mean = \rho q^2 + 2p$.
\end{itemize}
This together with the fourth equation in \eqref{four} leads to the choice of $p=\kappa_-$ in \eqref{p=kappa 2, 2d}.    With this choice, the identities $\mean = \rho q^2 + 2p$ and $\kappa = p\rho q^2+p^2$ become equivalent, so the fourth equation in \eqref{four} is then redundant.

\medskip
To proceed, we shall solve for $\rho$ and $v$ from the first three equations in \eqref{four}; or equivalently,
\begin{equation}\label{f eq, 2d}
\rho v_i v_j  = f_{ij}:= P_{ij} -\kappa_- g_{ij}.
\end{equation}
The crucial observation is that this system has further redundancies. Thus, \eqref{f eq, 2d} contains at most two independent equations for three unknowns $\rho, v_1, v_2$. Therefore, coupling  \eqref{f eq, 2d} with the conservation of mass~\eqref{E1} results in a determined (or under-determined) system.

\begin{lemma}\label{lem: redundancy}
The system~\eqref{f eq, 2d} has at most two independent equations for $(\rho, v_1, v_2)$.
\end{lemma}

\begin{proof}
We check $\det f = 0.$ Then $\rho v_1 v_2 = f_{12}$ can be deduced from $\rho(v_1)^2=f_{11}$ and $\rho (v_2)^2 = f_{22}$.

Indeed, the choice of $P$ in \eqref{P def, 2d} gives us
\begin{equation*}
f = \{f_{ij}\} = \begin{bmatrix}
g_{11} (H^2_2-\kappa_-) -  g_{12}H^1_2 && -g_{11}H^2_1 + g_{12}(H^1_1-\kappa_-)\\
g_{21}(H^2_2-\kappa_-) - g_{22} H^2_1 && -g_{21} H^2_1 + g_{22} (H^1_1 - \kappa_-)
\end{bmatrix}.
\end{equation*} 
A similar computation as for \eqref{kappa in 2d} leads to
\begin{align*}
\det \, f&= (\det\, g)\, \det (P'-\kappa_-I_2)\\
&= (\det\, g)\left\{\det \, H' - \kappa_-\,{\rm tr}\,H' + (\kappa_-)^2\right\}\\
&=(\det \,g)\left\{\kappa_- \kappa_+ - \kappa_-(\kappa_++\kappa_-)+(\kappa_-)^2\right\} = 0.
\end{align*}
\end{proof}

In particular, $\{f_{11}, f_{22}\}$ must have the same sign, which may be assumed positive without loss of generality. (If their signs are negative, one just reverses the orientation of $(M,g)$, which will result in interchanging the inward and outward unit normal vector fields and thus taking $(H,\kappa_-)$ to $(-H, -\kappa_+)$.)  Then, for $\rho>0$ we set
\begin{equation}\label{v_i, 2d formula}
v_i := \sqrt{\frac{f_{ii}}{\rho}}, \qquad i \in \{1,2\}.
\end{equation}
This clearly satisfies \eqref{f eq, 2d}. 

Finally, we specify the density $\rho>0$ via the continuity equation \eqref{E1}.  By virtue of \eqref{v_i, 2d formula} we have (no summation convention here) that
\begin{equation*}
0=\sum_{i=1}^2 \na_i  \sqrt{\rho f_{ii}} =\frac{1}{\sqrt{\det\,g}} \sum_{i=1}^2 \p_i \sqrt{(\det\,g) \rho f_{ii}}.
\end{equation*}
This can be solved \emph{locally} by the method of characteristics:
\begin{equation}\label{char, 2d, 1}
\rho\Big(z^1(t), z^2(t)\Big) = \rho\Big(z^1(0), z^2(0)\Big) \,\exp\bigg\{-2\int_0^t \bigg(\p_1\sqrt{(\det\,g) f_{11}}+\p_2\sqrt{(\det\,g) f_{22}}\bigg)\,\dd s\bigg\},
\end{equation}
where for each $i\in \{1,2\}$
\begin{equation}\label{char, 2d, 2}
\frac{\dd z^i}{\dd t} = \sqrt{(\det\,g) f_{ii}}
\end{equation}
as long as $\rho\left(z^1(0), z^2(0)\right)>0$. 

\begin{remark}
We emphasise the local nature of the solution formulae~\eqref{char, 2d, 1} and \eqref{char, 2d, 2}. It is in general not guaranteed that $\rho$ is defined globally on a non-simply-connected surface $(M,g)$. 
\end{remark}

\section{Two families of special solutions}\label{sec: examples}

The previous section provides a systematic way of producing fluid variables $(\rho, v, p)$ satisfying the steady Euler equation on $(M,g)$ from the second fundamental form $H$ of the isometric embedding $(M,g) \emb \R^3$. In this subsection, we describe two families of special solutions to the fluid PDE, found with the help of Lemma~\ref{lem: momentum + irrotational}.

\subsection{A stronger condition}
Our special solutions in this section shall be found by solving  the equation in Lemma~\ref{lem: momentum + irrotational}: 
\begin{equation}\label{L special sol}
\frac{\rho}{2}\mathcal{L}\left(\frac{\kappa_+ - \kappa_-}{\rho}\right) + \mathcal{L} \kappa_- =0
\end{equation}
for an \emph{arbitrary} first-order differential operator $\mathcal{L}$. In view of Lemma~\ref{lem: momentum + irrotational}, if the fluid variables satisfy \eqref{L special sol} together with the Bernoulli law $\rho q^2 = \kappa_+-\kappa_-$ and the conservation of mass, then the fluid is necessarily \emph{irrotational}, and the conservation of momentum is automatically satisfied. The requirement that \eqref{L special sol} holds for arbitrary $\mathcal{L}$ appears over-determined in general. Nevertheless, for certain surfaces $(M,g)$ with very special geometric properties, we may obtain explicit solutions for $\rho$.

Let us now take one step back: we first consider the PDE system consisting of irrotationality, conservation of mass~\eqref{E1}, and \eqref{new, flow speed q}, all reproduced below:
\begin{equation}\label{system}
\begin{cases}
\omega := \na_1 v_2 - \na_2 v_1 = 0,\\
\na_k(\rho v^k) = 0,\\
\rho = \hat{\rho}(q^2) := \frac{\kappa_+ - \kappa_-}{q^2}.
\end{cases}
\end{equation}
Here, $\kappa_\pm$ are viewed as given variables and $\left\{\rho, v^1, v^2\right\}$ as unknowns. We refer to the third equation as the \emph{Bernoulli law}, as it expresses density $\rho$ as a function of the flow speed squared, thus eliminating the pressure from the balance of momentum. As aforementioned, this system is weaker than \eqref{L special sol}, conservation of mass, and Bernoulli law put together.

\begin{remark}
As commented in Chen--Dafermos--Slemrod--Wang \cite[p.636]{cdsw}, on $\R^2$, the system~\eqref{system} is mathematically equivalent to the system of conservation of mass and momentum for smooth flows. The same holds on a surface $(M,g)$, as irrotationality also persists. Such equivalence breaks down in the presence of shocks, as vorticity may be created and mechanical energy may be converted to heat. The system~\eqref{system} is more popular among the aerodynamists due to its mathematical simplicity and its  analogy with  incompressible fluids. See also Bers \cite{bers}. 
\end{remark}

\subsection{Example 1: CMC surfaces}
If we set
\begin{equation}\label{cmc}
\rho = \kappa_+ + \kappa_- = \text{mean  curvature $\mean$} = \text{constant},
\end{equation}
then $\frac{\rho}{2}\mathcal{L}\left(\frac{\kappa_+ - \kappa_-}{\rho}\right) + \mathcal{L} \kappa_- = \frac{1}{2}\mathcal{L}(\kappa_+ + \kappa_-) = 0$. Conservation of momentum and irrotationality thus hold by Lemma~\ref{lem: momentum + irrotational}. Such $(M,g)$ is known as a \emph{CMC (constant mean curvature)} surface.

\begin{remark} If an oriented CMC surface $(M,g)$ has $\mean < 0$, then we should reverse the orientation to make the mean curvature positive, in view of the identification $\rho = \mean$ in \eqref{cmc}. 
\end{remark}

For the steady Euler equations~\eqref{E1} and \eqref{E2} corresponding to the CMC surface here, flow velocity $v$ is both divergence-free and curl-free; namely, $\dd$ and $\dd^*$-free on $(M,g)$. Thus, the differential 1-form $v^\sharp := v_i \dd x^i$ canonically dual to $v = v^j \p_j$ is harmonic, hence is smooth by elliptic regularity. By the aforementioned remarks in \cite[p.636]{cdsw}, the system~\eqref{system} is equivalent to the conservation laws \eqref{E1} and \eqref{E2}. We are in the case of \emph{incompressible} Euler equations.

By now we have not used the Bernoulli law, \emph{i.e.}, the third equation in \eqref{system}. Taking it into consideration, we will show that many CMC surfaces $(M,g)$ do not admit solutions to the steady Euler equations~\eqref{E1} and \eqref{E2}, under our specification of the fluid variables from the geometric variables. But, before embarking on ruling out candidates for CMC surfaces, we first present three simple affirmative examples. 

\begin{example}\label{example: S2, R2, cylinder}
Consider the case that both of the principal curvatures $\kappa_\pm$ are constant. Then we have three possible choices for $(M,g)$: the round sphere, the right cylinder, and the plane. 

\begin{enumerate}
\item
For $(M,g)=\R^2$, a harmonic $1$-form on $\R^2$ is constant. Also $\mean = \kappa_+ - \kappa_- =0$, so for any constant velocity $v$ the Bernoulli law is satisfied. Thus any constant vector field $v \in \G(T\R^2)$ is a Euler solution.

\item
For $(M,g) = r\mathbb{S}^2$, the sphere $\left\{ P \in \R^3:\,|P| = r\right\}$ for any $r>0$ equipped with the round metric, $v \equiv 0$ clearly satisfies \eqref{system}. In fact, in Proposition~\ref{propn: CMC} below we shall show that the zero solution is the only solution on the sphere.

\item
For $(M,g) = s \mathbb{S}^1 \times \R := \left\{\left(s\cos\theta,\,s\sin\theta,\,z\right)^\top \in \R^3:\,\theta \in [0,2\pi[,\,z\in\R\right\}$ for some $s>0$, the cylinder whose inclusion into $\R^3$ is an isometric embedding, one has $\kappa_+ = s^{-1}$ and $\kappa_-=0$. Then the Bernoulli law forces the flow speed $q$ to be $1$ for any $s>0$. There are four solutions: $v \in \{\pm \p_z, \,\pm\p_\theta\}$. Here $\p_z \equiv (0,0,1)^\top$ and $\p_\theta$ is the angular vectorfield given by $\p_\theta\big|_{(x,y,z)^\top} := \left(-\frac{y}{s},\,\frac{x}{s},\,0\right)^\top$ for any $(x,y,z)^\top \in s \mathbb{S}^1 \times \R \subset \R^3$. 
\end{enumerate}
\end{example}

Now we discuss the limitations on CMC surfaces imposed by the Bernoulli law together with irrotationality and mass conservation. We have the following
\begin{proposition}\label{propn: CMC}
Let $(M,g)$ be a smooth surface of constant mean curvature whose corresponding fluid variables satisfy \eqref{system}. Then $(M,g)$ cannot be any surface in the following list, unless we are in the simple cases given by Example~\ref{example: S2, R2, cylinder}:
\begin{enumerate}
\item
A minimal surface, \emph{i.e.}, $(M,g)$ on which $\mean \equiv 0$;
\item
Any simply-connected  closed surface;
\item
Any surface whose set of umbilical points $\mathcal{U}(M):=\left\{x \in M:\,\kappa_+(x)=\kappa_-(x)\right\}$ has an accumulation point.

\end{enumerate}
\end{proposition}

\begin{proof}
If (1) holds, then by the Bernoulli law $\kappa_+ \equiv \kappa_-$ on $(M,g)$, as the flow speed $q$ is finite. This together with $\mean \equiv 0$ yields  $\kappa_\pm \equiv 0$, hence $(M,g) = \R^2$, which is covered by Example~\ref{example: S2, R2, cylinder}. 

Next, suppose that (3) holds, \emph{i.e.}, the set of umbilical points $\mathcal{U}(M)$ is non-isolated. On this set we have $\rho = {0}/{q^2}$ by the Bernoulli law. If $q > 0$ on $\mathcal{U}(M)$, then $\rho \equiv 0$ on $M$ due to the CMC condition. This returns to the case~(1) above. If instead $q=0$ on $\mathcal{U}(M)$ (in which case we view the Bernoulli law on $\mathcal{U}(M)$ as undefined), then $v = 0$ on $\mathcal{U}(M)$. But since $v$ is harmonic thanks to irrotationality and conservation of mass,  $v$ is real-analytic on $(M,g)$. Then by analytic continuation $v \equiv 0$ on $M$. Then, as the mean curvature $\mean = \rho$ is finite, we must have $\kappa_+ \equiv \kappa_-$ by the Bernoulli law. So $(M,g)$ is totally umbilical, hence is $\R^2$ or $r\mathbb{S}^2$ for some $r>0$. Both cases have been dealt with in Example~\ref{example: S2, R2, cylinder}. 

Now consider (2), \emph{i.e.}, assume that $M$ is closed and simply-connected. As $v^\sharp$ is a harmonic $1$-form, it must be constant on $M$ by  de Rham cohomology. But the classification of closed surfaces shows that $M$ is a topological 2-sphere. Thus, by the hairy ball theorem, we see that $v \in \G(TM)$ must vanish somewhere. Thus $v$ is constantly zero. As before, for the Bernoulli law to make sense, we must have $\kappa_+ \equiv \kappa_-$ on $M$. This returns to Example~\ref{example: S2, R2, cylinder}. \end{proof}

Proposition~\ref{propn: CMC} rules out many candidates for CMC surfaces in our model from physical considerations. Here we do not intend to mean that these CMC surfaces are unnatural in any sense; we only suggest that the particular mapping we constructed which sends a solution $H$ to the Gauss--Codazzi equations~\eqref{gauss}--\eqref{codazzi} to the solution $(\rho, v, p)$ to the steady Euler system~\eqref{system} do not work for such CMC surfaces.

On the other hand, there is an abundance of higher-genus CMC surfaces which are not covered by Proposition~\ref{propn: CMC}; for example, the Wente torus in \cite{w} and the complete CMC surfaces obtained via PDE gluing by Kapouleas and others. See, for example, the recent survey by Breiner--Kapouleas--Kleene \cite{bkk}.

\emph{A priori} it is unclear if the above CMC surfaces admit solutions to \eqref{system}, as this system is overdetermined in general when $\rho = {\rm constant}$. Nevertheless, assuming existence, we may prove the following weak compactness result.

\begin{proposition}\label{propn: weak compactness}
Let $\left\{ \left(M_n, g_n\right)\right\}$ be a sequence of smooth closed CMC, non-minimal immersed surfaces in $\R^3$. Assume this sequence converges in the bi-Lipschitz sense to a limiting manifold which is possibly non-smooth. More precisely, there exist a manifold $M_\infty$ with a Lipschitz metric $g_\infty$ and bi-Lipschitz homeomorphisms $\psi_n: \left(M_\infty, g_\infty\right) \to \left(M_n, g_n\right)$ such that the bi-Lipschitz constants of $\psi_n$ tend to 1. Suppose that the principal curvatures $\kappa_{n, \pm}$ of $(M_n, g_n) \emb \R^3$ satisfy the integral bound:
\begin{equation}\label{new, integral bound}
\sup_{n \in \mathbb{N}}\int_{M_n}\left|\frac{\kappa_{n,+} - \kappa_{n,-}}{\kappa_{n,+} +  \kappa_{n,-}}\right|\,\dd V_{g_n} \leq K_0 < \infty
\end{equation}

In the above setting, consider a sequence of vector fields  $\{v_n\} \subset \G(TM_n)$ satisfying the system~\eqref{system} of irrotationality, conservation of mass, and Bernoulli law  subject to the identification~\eqref{cmc}. Then the corresponding vector fields $w_n := (\psi_n)^\# v_n \in {\bf Lip}(M_\infty; TM_\infty)$ on the limiting manifold converge in the weak $L^2$-topology to some $\overline{w} \in L^2(M_\infty; TM_\infty)$. 

Here, the weak limit $\overline{w}$ is divergence-free and irrotational in the sense of distributions. In addition, $\left|\overline{w}\right|_{g_\infty}^2$ is the distributional weak limit of $(q_n)^2 := \left|v_n\right|_{g_n}^2$ on $\left(M_\infty, g_\infty\right)$.   
\end{proposition}

Before giving a proof, we first remark on the statement of Proposition~\ref{propn: weak compactness}. Therein, the sequence of CMC surfaces $\left\{ \left(M_n, g_n\right)\right\}$ may be obtained, for example, from scaling the metrics for a given CMC surface. The limiting manifold $\left(M_\infty, g_\infty\right)$ is not required to be a Riemannian manifold, but some regularity is assumed. Indeed, it should at least be a Lipschitz manifold with Lipschitz metric; in particular, it is assumed to be non-collapsing. Because $\psi_n$ are bi-Lipschitz homeomorphisms, the pullback $w_n := (\psi_n)^\# v_n$ are well defined Lipschitz vector fields on $M_\infty$. One subtlety is that we have assumed $(M_n,g_n)$ to be smooth manifolds (intrinsically), but their CMC immersions in $\R^3$ may fail to be smooth. This is the case for the Wente torus in \cite{w}.

The key point of the theorem is to establish the weak $L^2$-convergence  $w_n \weak \overline{w}$, such that the limiting vector field $\overline{w}$ also satisfies the system~\eqref{system} in a ``very weak sense''.  By this we mean that the irrotationality and conservation of mass hold weakly. Moreover, although the principal curvatures $\kappa_{n,\pm}$ of each CMC surface $(M_n, g_n) \emb \R^3$ may blow up as $n \nearrow \infty$, the Bernoulli law and  \eqref{new, integral bound} implies that $\left\{(q_n)^2 = \frac{\kappa_{n,+} - \kappa_{n,-}}{\kappa_{n,+} +  \kappa_{n,-}}\right\}$ is uniformly bounded in $L^1$. Note that for the integral bound condition~\eqref{new, integral bound} to hold, $(M_n, g_n)$ are necessarily non-minimal. 


\begin{proof}
Let us work with differential 1-forms instead of vector fields. Let $\alpha_n := (v_n)^\sharp$ be the 1-forms canonically dual to $v_n$ via metrics $g_n$. In this case, the irrotationality and the conservation of mass of $v$ are equivalent to $\dd \alpha_n=0$ and $\dd^* \alpha_n =0$. Equivalently, $\alpha_n$ are harmonic 1-forms on the smooth Riemannian manifolds $(M_n, g_n)$. Note that they are smooth and are confined in finite-dimensional spaces (dimension = the first betti number of $M$).

Consider the Lipschitz 1-forms:
\begin{align*}
\beta_n := (\psi_n)^\# \alpha_n = (v_n)^\sharp.
\end{align*}
We \emph{claim} that $V_0:=\sup_{n \in \mathbb{N}}\|\beta_n\|_{L^2(M_\infty,g_\infty)}$ is finite. Indeed, by the duality $TM_n \cong T^*M_n$, Bernoulli law, and the integral bound~\eqref{new, integral bound}, we have
\begin{align*}
\|\alpha_n\|_{L^2(M_n,g_n)} := \|q_n\|_{L^2(M_n,g_n)} = \sqrt{\left\|\frac{\kappa_{n,+} - \kappa_{n,-}}{\kappa_{n,+} +  \kappa_{n,-}}\right\|_{L^1(M_n,g_n)}\,} \leq \sqrt{K_0}.
\end{align*}
So $V_0<\infty$ depends only on $K_0$ and the uniform bound on bi-Lipschitz constants of $\psi_n$.  As $\dd \alpha_n$ and the differential commutes with pullbacks, we have that
\begin{equation}\label{d beta = 0}
\dd \beta_n = 0\qquad \text{ on } \left(M_\infty, g_\infty\right).
\end{equation}

Now let us turn to the divergence of $v_n$; that is, the co-differential of $\beta_n$. We have
\begin{align*}
\dd^* \beta_n = \star \dd \star \left( (\psi_n)^\# \alpha_n\right),
\end{align*}
where $\star$ is the Hodge star. Considering the  commutator $T_n := \left[\star \dd\star,\,(\psi_n)^\# \right]$, we have
\begin{align*}
\dd^* \beta_n = T_n \alpha_n + (\psi_n)^\# \left(\star \dd \star  \alpha_n\right) = T_n \alpha_n,
\end{align*}
since $\dd^*\alpha_n = 0$. But $T_n$ is a pseudo-differential operator of order zero with $L^\infty$-coefficients, so 
\begin{equation}\label{d star beta = 0}
\left\{\dd^*\beta_n\right\} \quad \text{ is bounded in $L^2\left(M_\infty, g_\infty\right)$}.
\end{equation}

By a standard application of the div-curl lemma (see Murat \cite{m1, m2} and Tartar \cite{t1, t2}, or its adaptations on manifolds in \cite{cl, clnew}), \eqref{d beta = 0} and \eqref{d star beta = 0} implies that $\beta_n \weak \overline{\beta}$ in the weak $L^2$-topology to a limiting 1-form $\overline{\beta} \in L^2(M_\infty; T^*M_\infty)$, such that
\begin{align*}
\bra \beta_n, \beta_n \ket_{g_\infty} \, \longrightarrow \, \bra \overline{\beta},\overline{\beta}\ket_{g_\infty} \text{ in the sense of distributions}. 
\end{align*}
Therefore, by setting $\overline{w} := \left(\overline{\beta}\right)^\flat \in  L^2(M_\infty; TM_\infty)$ and invoking the assumption that the bi-Lipschitz constants of $\psi_n$ converge to 1, we conclude that 
\begin{align*}
q_n^2 \, \longrightarrow \, \left|\overline{w}\right|^2_{g_\infty} \text{ in the sense of distributions}.
\end{align*}
This completes the proof.   \end{proof}

\begin{remark}
The Wente torus ${\bf T}_W$, \emph{cf.} \cite{w}, serves as an illustrative example for Proposition~\ref{propn: weak compactness}: its metric is $g=e^{2\omega}\delta$ (where $\delta$ is the Euclidean metric) and its principal curvatures are $\kappa_+ = e^{-\omega} \cosh \omega$ and $\kappa_- = e^{-\omega}\sinh\omega$, hence $\mean =1$ and $q^2 = e^{-2\omega}$. 

Consider the  homothety $\iota \rightsquigarrow \eta\iota$ for some $\eta>0$: here $\iota$ is the isometric immersion of ${\bf T}_W$ into $\R^3$. Then the metric associated to $\eta\iota$ is 
$g_\eta := \eta e^{2\omega}\delta$, and its second fundamental form is $\eta H_{ij}$, where $H_{ij}$ is the second fundamental form associated to $\iota$. We are interested in sending $\{\eta_n\} \to \eta_0$, where $\eta_n$ and $\eta_0$ are positive constants. Along this limiting process, the principal curvatures remain unchanged, hence the flow speeds $q_{\eta_n}$ for $\left({\bf T}_W, g_{\eta_n}\right)$ are the same for all $n$, thanks to the Bernoulli law in \eqref{system}. Then, assuming the existence of solution on $\left({\bf T}_W, g_{\eta_n}\right)$ for each $\eta_n$, we find that the hypotheses in Proposition~\ref{propn: weak compactness} are immediately verified. Indeed, for \eqref{new, integral bound}, the integral on the left-hand side equals $\left\|q_{\eta_n}\right\|_{L^2\left({\bf T}_W, g_{\eta_n}\right)} = \left(\eta_n\right)^2 \int_{{\bf T}_W} e^{-2\omega}\,\dd V_g$, where $\omega$ is a smooth function on the compact Wente torus $({\bf T}_W,g)$, hence bounded.  
\end{remark}

\subsection{Example~2: Standard torus}\label{subsec: torus}

In this example, we consider a 2-dimensional connected, closed, orientable manifold $(M,g)$ and look for solutions with $\kappa_+ = {\rm constant}$ and $\rho = \widetilde{\rho}(\kappa_-)$. By purely geometrical arguments, such surfaces can be completely classified: Shiohama--Takagi \cite{st} proved that except for the three surfaces listed in Example~\ref{example: S2, R2, cylinder}, the only possible such $(M,g)$ is the \emph{standard torus}; namely, the surface of revolution
\begin{equation}\label{torus}
\tor = \left\{\begin{bmatrix}
(a+b\cos \theta)\cos \phi\\
(a+b\cos \theta)\sin \phi\\
b\sin\theta
\end{bmatrix}:\,0 \leq \theta,\phi < 2\pi\right\}.
\end{equation}
Here $a$ and $b$ with $a>b>0$ are the major and minor radii of $\tor$. The metric on $\tor$ is the pullback of the Euclidean metric on $\R^3$ by the inclusion $\tor \subset \R^3$, denoted by $g$ in this example. 

Note that the standard torus has non-constant mean curvature (and, in addition, sign-changing Gauss curvature), so the corresponding fluid satisfying the steady Euler equations cannot be  incompressible. Conversely, if $\rho$ is non-constant, then $\kappa_-$ is non-constant, so we are in neither of the three cases in Example~\ref{example: S2, R2, cylinder}. Thus the classification result in \cite{st} shows that $(M,g)=\tor$ for some $a>b>0$.

When $\rho = \widetilde{\rho}(\kappa_-)$, the equation in Lemma~\ref{lem: momentum + irrotational} becomes
\begin{align*}
\frac{1}{2}\left\{1-\frac{\kappa_+-\kappa_-}{\rho}\cdot \widetilde{\rho}\,'(\kappa_-) \right\} \mathcal{L} \kappa_- = 0.
\end{align*}
Then it suffices to solve the ODE in the $\kappa_-$ variable:
\begin{align*}
1-\frac{\kappa_+-\kappa_-}{\widetilde{\rho}(\kappa_-)}\cdot \widetilde{\rho}\,'(\kappa_-) = 0.
\end{align*}
To this end, we make the change of variable
\begin{equation*}
F(\delta) := \log \widetilde{\rho}(\kappa_-),\qquad \delta:=\kappa_+-\kappa_- > 0.
\end{equation*}
In the case that $\kappa_-$ is non-constant, we have
\begin{align*}
\kappa F'(\kappa) = -1,
\end{align*}
thus $F(z) = C - \log z$ for an arbitrary constant $C$. Equivalently, the solution is
\begin{equation}\label{rho solution}
\widetilde{\rho}(\kappa_-) = \frac{C_0}{\kappa_+-\kappa_-},
\end{equation}
where $C_0>0$ is an arbitrary constant, $\kappa_+ > \kappa_-$, and $\kappa_+ = {\rm constant}$.

Now let us analyse the explicit solution \eqref{rho solution} from the perspective of fluid dynamics. Recall from Step~2, \eqref{p=kappa 2, 2d} that $\kappa_- = p$. So the density $\rho$ and pressure $p$ are related by $\rho = \frac{C_1}{C_2-p}$ for constants $C_1>0$ and $C_2>p$, or alternatively
\begin{equation}\label{chaplygin}
p = C_2 - \frac{C_1}{\rho}\quad\text{ and } \quad p'(\rho) = \frac{C_1}{\rho^2} =:c^2\quad \text{ ($c$ is the sonic speed)}.
\end{equation}
Since we wish $p$ to be increasing in $\rho$, we choose $C_1>0$. This resembles the constitutive relation for the \emph{Chaplygin gas}; see Chen--Slemrod--Wang \cite[p.415]{csw}. In this case, the steady fluid motion is \emph{sonic}: one infers from the Bernoulli law (\emph{i.e.}, the third equation in \eqref{system}) that 
\begin{align}\label{torus sonic}
q^2 = \frac{\kappa_+-\kappa_-}{\rho} = \frac{C_0}{\rho^2} = \frac{d}{d\rho} \left(\kappa_+ - \frac{C_0}{\rho}\right) = p'(\rho),
\end{align}
hence the Mach number $M=q/c = 1$.  Note that the relation \eqref{torus sonic} agrees with the classical Bernoulli law for steady flow of polytropic gases (Courant--Friedrichs \cite[p.22, (14.05)]{cf}): $$q^2 +\frac{2}{\gamma-1}c^2= \widehat{q}^2.$$ That is, $$q^2 - p'(\rho) = \widehat{q}^2$$ when applied to the Chaplygin gas, and the Bernoulli constant $\widehat{q}$ is taken to be zero. This of course yields $$q^2\rho^2 = C_1.$$

We now substitute the Bernoulli law  into the conservation of mass~\eqref{E1}, in provision that the irrotationality condition is already satisfied. The system~\eqref{system} is invariant under the scaling
\begin{equation}\label{scaling}
(\rho, v) \longmapsto \left( \lambda^2 \rho,\,\frac{1}{\lambda}v \right) \qquad \text{ for each constant } \lambda > 0.
\end{equation}
Thus, without loss of generality, we may fix $C_0=1$ in \eqref{rho solution} or \eqref{torus sonic}, which leads to 
\begin{align*}
\rho = \frac{1}{\kappa_+-\kappa_-}\quad\text{and}\quad |v|_g = q = \kappa_+-\kappa_-.
\end{align*}
As a consequence, $\na_k\left(\rho v^k\right) = 0$ becomes the  equation:
\begin{equation}\label{1-laplace}
{\rm div}_g\left( \frac{v}{|v|_g} \right) = 0,
\end{equation}
where ${\rm div}_g: \G(TM) \to C^\infty(M)$, $u \mapsto \na_j u^j$ is the Riemannian divergence with respect to $g$. 

In summary, we have shown the following
\begin{lemma}\label{lem: system standard torus}
On the standard torus $(M,g) = \tor$ for $a>b>0$, for the ansatz $\rho = \widetilde{\rho}(\kappa_-)$, the system~\eqref{system} of irrotationality, conservation of mass, and Bernoulli law  is equivalent, modulo a multiplicative constant due to the scaling invariance in \eqref{scaling}, to the following system:
\begin{equation*}
\begin{cases}
{\rm div}_g\left(\frac{v}{q}\right) = 0,\\ 
q = \kappa_+ - \kappa_-, \\ \rho q =1.
\end{cases}
\end{equation*}
\end{lemma}

\begin{remark}\label{rem: no stagnation}
In this lemma, thanks to the Bernoulli law $q^2 = \frac{\kappa_+ - \kappa_-}{\rho}$, the flow has no stagnation point (\emph{i.e.}, $q \neq 0$ everywhere), for otherwise there would be an umbilical point on $(M,g)$, which would lead to $(M,g) = \R^2$ or $r\mathbb{S}^2$ but not $\tor$ in the presence of one constant principal curvature. See Shiohama--Takagi \cite{st}. 
\end{remark}

Special solutions with  geometric significance can be found for the system in Lemma~\ref{lem: system standard torus}. Indeed, as computed in  Appendix~\ref{sec: app A}, by expressing $v \in \G(T\tor)$ as
\begin{align*}
v = v^\theta (\theta,\phi) \p_\theta + v^\phi (\theta,\phi) \p_\phi
\end{align*}
in the natural co-ordinate system $\{\p_\theta, \p_\phi\}$ for  the parametrisation of $\tor$ in \eqref{torus}, we may write the first-order system in Lemma~\ref{lem: system standard torus} as follows:
\begin{equation}\label{local, 1-Laplace system}
\begin{cases}
-{2b}\sin{\theta} \cdot v^\theta(\theta,\phi) + \left( a+ b\cos {\theta} \right) \left[\p_\theta v^\theta(\theta,\phi) + \p_\phi v^\phi(\theta,\phi)\right] = 0,\\
\left( a+ b\cos {\theta} \right)^2 \left(v^\theta(\theta,\phi)\right)^2 + {b^2}\left(v^\phi(\theta,\phi)\right)^2 = \frac{a^2}{\left( a+ b\cos {\theta} \right)^2}.
\end{cases}
\end{equation}

\begin{enumerate}
\item
By inspection, the ansatz
\begin{equation*}
v^\theta \equiv 0,\qquad v^\phi(\theta,\phi) = \widehat{v^\phi}(\theta)
\end{equation*}
verifies the first equation in \eqref{local, 1-Laplace system}. Then we can directly determine $\widehat{v^\phi}(\theta)$ from the second equation. In this way, we obtain a solution:
\begin{equation}\label{special sol a, torus}
v = \pm \frac{a}{b^2\left(a+b\cos\theta\right)^2} \p_\phi.
\end{equation}
Notice that \eqref{special sol a, torus} yields a well defined vector field $v \in \G(T\tor)$, since $\p_\phi$ is globally defined over $\tor$ and this expression is $2\pi$-periodic in $\theta$. Integral curves of $v$ are precisely the \emph{toroidal geodesics}. The magnitude of $v$ is constant on each of the toroidal geodesics but varies in the poloidal direction.
\item
On the other hand, the ansatz
\begin{align*}
v^\phi \equiv 0,\qquad v^\theta(\theta,\phi) = \widehat{v^\theta}(\theta)
\end{align*}
also yields a solution. The profile $\widehat{v^\theta}$ is determined by the second equation of \eqref{local, 1-Laplace system}, and one may directly check that it verifies the first equation. The solution is
\begin{equation}\label{special sol b, torus}
v = \pm \frac{a}{b(a+b\cos\theta)^2}\,\p_\theta.
\end{equation} 
It is rotationally symmetric; \emph{i.e.}, $v$ is completely determined by its profile on any of the generating circles $\{\phi = {\rm const.}\}$ of the $\tor$ as a surface of revolution. The integral curves of $v$ are precisely the poloidal geodesics. 
\end{enumerate}
In both cases (1) $\&$  (2) above, the density is given by
\begin{equation}
\rho = \rho(\theta) =  \frac{b\left( a+b\cos\theta \right)^2}{a}. 
\end{equation}

\begin{remark}
We emphasise that the topology of $\tor$ plays a crucial role in the discovery of the special solutions~\eqref{special sol a, torus} and \eqref{special sol b, torus}, although in a somewhat fortuitous manner. Indeed, these expressions rely on the properties that the standard torus admits everywhere non-vanishing vector fields, and that the toroidal and poloidal geodesics respectively foliates the whole surface. 

In contrast, such properties do not hold on topological 2-sphere $\Sigma$ due to the hairy ball theorem. The Chaplygin-type Bernoulli law $$\rho q = {\rm constant}$$ cannot hold everywhere  on $\Sigma$, as stagnation points are unavoidable. Thus, \emph{even without assuming the stronger condition~\eqref{L special sol}} and only considering the system~\eqref{system} of irrotationality, conservation of mass, and Bernoulli law, it is perhaps undesirable to impose a Chaplygin-type Bernoulli law.  

\end{remark}

\section{Irrotational Chaplygin gas on the flat torus}\label{sec: Irrotational Chaplygin gas on the flat torus}

In the above section \S\ref{subsec: torus}, we investigated the irrotational, steady Euler equation for the Chaplyngin-type gas (\emph{i.e.}, with the Bernoulli law $\rho q = {\rm constant}$) on the standard torus $\tor$. Now we discuss the analogous issue on the  \emph{flat} torus $$\ttwo = [0,1] \times [0,1]\slash \sim.$$ Here $\sim$ is the equivalence relation on $[0,1]\times [0,1]$ given by $(0,y) \sim (1,y)$ for all $y \in [0,1]$ and $(x,0) \sim (x,1)$ for all $x \in [0,1]$. The quotient space inherits the Euclidean metric from $[0,1]\times [0,1]$ and is a Riemannian manifold.

Our hope is to provide some new insights into the irregular isometric embeddings of $\ttwo$ into $\R^3$ constructed by Nash \cite{nash}, Kuiper \cite{k1, k2}, Borisov \cite{bor}, Borrelli--Jabrane--Lazarus--Thibert \cite{b1, b2}, and Conti--De Lellis--Sz\'{e}kelyhidi Jr. \cite{cds}, among many others. In particular, it is well-known (see the above references) that such irregular isometric embeddings are closely related to the irregular or ``wild'' solutions to the Euler equations in fluid dynamics, especially in terms of the methods of their constructions. It is thus natural to study the Euler equations on the flat torus $\ttwo$ and discuss its possible implications on the irregular isometric embeddings.

\subsection{From standard torus to flat torus: convergence of metrics} \label{subsec: From standard torus to flat torus: convergence of metrics}

Our motivating observation for the developments in this section is as follows. Consider a rescaled version of the standard torus $\tor$ (note that $\mathbf{T}_1(a,b)=\tor$): 
\begin{align}\label{torc}
\torc  = \left\{\begin{bmatrix}
\left[a+b\cos \left(\frac{\theta}{c}\right)\right]\cos \phi\\
\left[a+b\cos \left(\frac{\theta}{c}\right)\right]\sin \phi\\
b\sin\left(\frac{\theta}{c}\right)
\end{bmatrix}:\,0 \leq \phi < 2\pi,\,0 \leq \theta < 2\pi c\right\}\quad \text{where } c>0.
\end{align}
We are interested in the limiting process $b, c \searrow 0$. A straightforward computation yields that the metric $g = g_{a,b,c}$ for $\torc$ is given by
\begin{equation}\label{metric, torc, maintext}
g = \begin{bmatrix}
\left[a+b \cos \left( \frac{\theta}{c} \right)\right]^2 & 0\\
0&\frac{b^2}{c^2}
\end{bmatrix}.
\end{equation}

\begin{remark}
The rescaled tori $\torc$ are motivated by the ``Nash wrinkles'', which are the basic building block for the $C^1$-isometric embeddings  $\ttwo \emb \R^3$. Indeed, the Nash wrinkles contain terms of the form ${\e} \cos\left(\frac{\theta}{\e}\right)$ or ${\e} \sin\left(\frac{\theta}{\e}\right)$ in the approximation sequence  to isometric embeddings for $0 < \e \ll 1$. See \cite[Equations~(13) $\&$ (15)]{nash}.
\end{remark}

 If we fix ${b}/{c} = {\rm constant} = c_0$ and send  $b, c \searrow 0$, then:
\begin{enumerate}
\item
 $g$ tends to the flat metric $\overline{g} = \begin{bmatrix}
a^2 & 0 \\
0 & (c_0)^2
\end{bmatrix}$ in $L^\infty$;
\item
the first derivative $\p g$ converging in the weak-$\star$ topology of $L^\infty$ to zero. 
\end{enumerate}
The second property holds since $$\p_\theta g_{11} = -2\left[a+b \cos \left( \frac{\theta}{c} \right)\right]\sin \left( \frac{\theta}{c} \right),$$ where $b \searrow 0$ and $\sin(\theta\slash c) \weak 0$ weakly-$\star$ in $L^\infty$ by the Riemann--Lebesgue lemma as $c \searrow 0$. Therefore, we expect the behaviour of the solutions to the Euler equations (in the form of \eqref{system}) along this limiting process to contain certain persistent information about the steady irrotational flows, which in turn will single out ``good'' flows on the flat torus.

To this end, we carry out the process in \S\ref{subsec: torus} for the rescaled standard torus $\torc$. The computations are exactly parallel (see Appendix~\ref{sec: app A}). The principal curvatures of $\torc$ are
\begin{equation*}
\kappa_+ = \frac{1}{b}\quad \text{ and } \quad \kappa_- = \frac{\cos  \left( \frac{\theta}{c} \right)}{\left[a+b \cos \left( \frac{\theta}{c} \right)\right]},
\end{equation*}
and as before (see \eqref{p=kappa 2, 2d}) we take $p = \kappa_-$. By seeking solutions of the form $\rho = \widetilde{\rho}(\kappa_-)$ we again arrive at the Bernoulli law of Chaplygin-type gas. Therefore, we once more obtain the PDE system in Lemma~\ref{lem: system standard torus}, reproduced below:\begin{equation*}
\begin{cases}
{\rm div}_g\left(\frac{v}{q}\right) = 0,\\ 
q = \kappa_+ - \kappa_-, \\ \rho q =1.
\end{cases}
\end{equation*}
Thus
\begin{align}\label{special sol on Torc}
q= \frac{a}{b\left[a+b \cos \left( \frac{\theta}{c} \right)\right]^2}\quad\text{ and } \quad \rho = \frac{b\left[a+b \cos \left( \frac{\theta}{c} \right)\right]^2}{a}.
\end{align}
Moreover, rewriting the above system in local co-ordinates (see \eqref{local, 1-Laplace system} and ensuing developments for $\tor = \mathbf{T}_1(a,b)$), we obtain  special solutions analogous to  those in \eqref{special sol a, torus} and \eqref{special sol b, torus}:
\begin{equation}\label{special solutions, torc, v}
v=\pm \frac{ac}{b^2\left[a+b\cos \left( \frac{\theta}{c}\right)\right]} \p_\phi\quad\text{ or } \quad\pm \frac{a}{b\left[a+b\cos \left( \frac{\theta}{c}\right)\right]^2}\,\p_\theta.
\end{equation}

In the limit $b, c \searrow 0$ while keeping ${b}/{c} = {\rm constant}$, we see that $q \to \infty$ and $\rho \to 0$, both in the pointwise sense. In terms of geometric variables, we have $\kappa_+ \nearrow +\infty$ and $\kappa_- \weak 0$ in the weak-$\star$ topology of $L^\infty$, thanks to the presence of rapid oscillations. Moreover, in view of the identification between the stress-momentum tensor $P$ and the second fundamental form $H$ (see \eqref{P def, 2d}, $H \approx \rho v\otimes v + g\kappa_-$), the extrinsic geometry of the associated isometric embeddings blows up pointwise in the limit. Nevertheless, the relation $\rho q =1$ remains to hold everywhere.

The above investigations hint at the following:
\begin{quote}
On the flat torus $\ttwo = [0,1]\times [0,1] \slash \sim$, it is natural to consider the steady, irrotational Euler equations for the Chaplygin gas, namely that
\begin{equation}\label{flat torus eq}
\begin{cases}
\omega = 0,\\
{\rm div}_g(\rho v) = 0,\\
\rho q =1.
\end{cases}
\end{equation}
\end{quote}

See \S\ref{subsec: geom remark} below for  further remarks. 

\subsection{On the flat torus}
Now we turn to the study of \eqref{flat torus eq} on $\ttwo$.    Substituting the Bernoulli law $\rho q =1$ into the conservation of mass, we obtain
\begin{equation}
{\rm div}_g \left(\frac{v}{q}\right) = 0\qquad \text{ on } \ttwo.
\end{equation}
In local co-ordinates the Riemannian divergence is given by 
\begin{align*}
{\rm div}_g w = \frac{1}{\sqrt{\det g}} \p_i \left( \sqrt{\det g} w^i \right)
\end{align*}
for any vector field $w$. But $\det g$ is constant for the flat metric, so the equation becomes
\begin{align*}
\p_1 \left( \frac{v^1}{|v|_g}\right) + \p_2 \left( \frac{v^2}{|v|_g}\right) = 0.
\end{align*}
Clearly, we have the following special solutions which are also irrotational:
\begin{equation}\label{shear flow special sol}
v = f\left(x^1\right)\p_1 \quad \text{ or } \quad v = h\left(x^2\right)\p_2,
\end{equation}
where $f, h \in C^\infty([0,1])$ are arbitrary 1-periodic smooth functions without zeros. These are well defined global vector fields on $\ttwo$, and they represent horizontal or vertical flows with arbitrary smooth, stagnation point-free velocity profiles. Throughout we  identify $\ttwo = [0,1] \times [0,1] \slash \sim$, with $\p_1$, $\p_2$ tangent to each copy of $[0,1]$.

In view of the irregular isometric embedding $\iota: \ttwo \emb \R^3$ remarked at the beginning of this section (\emph{cf}.  \cite{nash, k1, k2, bor, b1, b2, cds}), we  obtain a simple recip\'{e} for irregular steady Euler solutions which exhibit fractal patterns of motion. One may take $\iota \in C^{1,\beta}(\ttwo,\R^3)$; the best $\beta$ known to date is $\frac{1}{7}-\e$. See \cite{bor, cds}.

\begin{theorem}\label{thm: fractal euler}
Let ${\bf T}$ be the embedded $C^{1,\beta}$-surface in $\R^3$ which is an isometric copy of the flat torus $\ttwo$. There exists a $C^{1,\beta}$-vector field $V \in \G(T\R^3) $ whose restriction to ${\bf T}$ satisfies the Euler equation~\eqref{flat torus eq} for steady, irrotational Chaplygin gas.
\end{theorem}

\begin{proof} 
Let $v$ be the smooth shear flows defined in \eqref{shear flow special sol} on $\ttwo$, and let $\iota: \ttwo \emb \R^3$ be an  $C^{1,\beta}$-isometric embedding. The pushforward $\iota_\ast v$ defines a solution on $\iota({\bf T})$, which can be extended to a   $C^{1,\beta}$-vector field $V$ by Whitney embedding.   \end{proof}

\begin{remark}
The velocity field $v \in \G(T\ttwo)$ is ``intrinsically smooth'', in the sense that it is a smooth vector field on $\ttwo$, which is itself a smooth Riemannian manifold. The irregularity arises only from the extrinsic geometry; namely, the isometric embedding $\iota: \ttwo \emb \R^3$.  
\end{remark}

\subsection{Geometric remarks}\label{subsec: geom remark}

Earlier in this section we presented some investigations on the irregular solutions to the steady motion of Chaplyngin gas on the flat torus. Such discussions are of a  tentative and preliminary nature; we shall collect here a few related remarks. 

Our first observation is that the smooth solutions $v \in \G(T\ttwo)$ in \eqref{shear flow special sol}, or the $C^{1,\beta}$-solutions in Theorem~\ref{thm: fractal euler}, cannot arise from the limiting process introduced in \S\ref{subsec: From standard torus to flat torus: convergence of metrics}. Indeed, as computed in \eqref{special sol on Torc}, the flow speed becomes unbounded as $b, c \searrow 0$. This is because $q = \kappa_+ - \kappa_-$, where the larger principal curvature $\kappa_+ = 1/b$ blows up, while the smaller principal curvature is oscillatory (weakly-$\star$ converges to zero in $L^\infty$).

For the limiting process in \S\ref{subsec: From standard torus to flat torus: convergence of metrics}, we have shown that the metrics $g=g_{a,b,c}$ for $\torc$ converge in the $L^\infty$-norm to the flat metric $\overline{g}$ on $\ttwo$ (modulo scaling by constants). This, however, does not imply that the sequence of manifolds $\left\{\torc\right\}$ converges to $\ttwo$ in the Hausdorff--Gromov sense. Indeed,  $\torc$ converge to the central toroidal circle
\begin{align*}
{\bf C}_a = \left\{\begin{bmatrix}
a \cos\phi\\
a \sin\phi\\
0
\end{bmatrix}:\, 0 \leq \phi < 2\pi\right\}
\end{align*}
as $b \searrow 0$. (The limit $c \searrow 0$ has the effect that one traverses the poloidal circles, \emph{i.e.}, in the $\theta$-direction, increasingly rapidly.)  One possible heuristic    for this phenomenon is  to view $\mathbf{C}_a$ as a ``degenerate flat torus'', namely that $\ttwo \cong {\bf C}_a \times \mathbb{S}^1$ with the $\mathbb{S}^1$ factor collapsing to a point. For $0 < b \approx c \ll 1$ we may view the two special solutions found in \eqref{special solutions, torc, v} as being concentrated  on the central toroidal circle ${\mathbf C}_a$ and on the vanishing poloidal direction $\p_\theta$, respectively. This is vaguely reminiscent of the 5-dimensional Kaluza--Klein theory of cosmology (see \emph{e.g.},  \cite{5d}), which has an extra dimension curled up into a small loop in addition to the 4-dimensional spacetime.

In view of the previous paragraph, it  is natural to consider other approximating sequences of metrics to the flat torus, as we are not enforcing  Gromov--Hausdorff convergence of manifolds. For example, one may take on an arbitrary surface $M$ the metrics (as long as they exist on $M$):
\begin{align*}
g_n = \begin{bmatrix}
1 + \e_1(n) & 0\\
0 & 1 + \e_2(n)
\end{bmatrix}, \qquad \e_1(n), \e_2(n) \longrightarrow 0 \text{ in suitable sense as } n \to \infty.
\end{align*}

On an arbitrary closed (\emph{i.e.}, compact, without boundary) manifold of arbitrary dimensions, one can always solve (in the distributional sense) the system of irrotationality and  conservation of mass for Chaplygin gas, away from the set of stagnation points. If there exist zero-free  harmonic 1-forms, then one can further require such solutions to be zero-free. 
This is summarised as the following theorem, whose proof is presented in Appendix~\ref{sec: appendix, variation}.

\begin{theorem}\label{thm: chaplygin}
Let $(M,g)$ be any closed Riemannian manifold. Let $h \in \Omega^1(M)$ be any nontrivial harmonic 1-form. There exist $\rho: M \to [0,\infty[$ and $v \in L^1(M;TM)$ which are weak solutions to the steady Euler equations away from the set of stagnation points of $v$: 
\begin{equation*}
\begin{cases}
\omega = 0,\\
{\rm div}_g(\rho v) = 0,\\
\rho q = {\rm constant},
\end{cases}
\end{equation*}
such that $v^\sharp$, the 1-form canonically dual to $v$, lies in the same cohomology class as $h$. If, in addition, $(M,g)$ admits an everywhere non-vanishing harmonic 1-form, then $v$ can be chosen without stagnation points.
\end{theorem}

\section{From Nash--Kuiper to Euler: 2D Irregular Solutions and Renormalisation}\label{sec: renorm}

We have described a connection between the smooth, steady compressible Euler equations and the isometric embedding of surfaces. Such constructions cannot be directly applied to irregular (non-$C^2$) embeddings obtained from the Nash--Kuiper Theorem~\ref{thm: nash kuiper}, since neither the second fundamental form nor the Gauss curvature are well defined in the little regularity setting.

In this section, we establish a link between irregular isometric embeddings and compressible Euler equations. We bypass the aforementioned hindrance by a ``renormalisation'' process, based on explicit estimates in the corrugation/convex integration process of Conti--De Lellis--Sz\'{e}kelyhidi Jr. \cite{cds}. Indeed, in \cite{cds} the authors  provided a constructive proof for a version of the Nash--Kuiper theorem, which yields the best H\"{o}lder exponent $\alpha$ up to date. In what follows, unless otherwise specified, $\|\bullet\|_k$ for $k=0,1,2,\ldots,\infty$ denotes the $C^k$-norm.
\begin{theorem}[Theorem 2 in \cite{cds}]\label{thm: cds}
Let $M$ be an $n$-dimensional
compact manifold with a Riemannian metric $g$ in $C^\beta(M)$ and $m\geq n+1$. Then there is a constant $\delta_0 > 0$ such that if $u \in C^2(M;\R^{m})$ and $\alpha$ satisfy
\begin{equation*}
\|\p_i u \cdot \p_j u -g_{ij}\|_0 \leq \delta_0^2\qquad \text{ and } \qquad 0<\alpha < \min \bigg\{ \frac{1}{2(n+1)n_\star},\frac{\beta}{2}\bigg\}
\end{equation*}
where $n_\star = n(n + 1)/2$,  then there exists a map $v \in C^{1,\alpha}(M;\R^m)$ such that
\begin{equation*}
\p_i v \cdot \p_j v =g_{ij} \qquad \text{ and } \qquad \|v-u\|_1 \leq {\rm const.} \max_{i,j} \sqrt{\|\p_i u \cdot \p_j u -g_{ij} \|_0}.
\end{equation*} 
\end{theorem}
 The proof of  Theorem \ref{thm: cds} involves the construction of a sequence of smooth maps $\{u_q\}_{q\in\mathbb{N}}$ via the ``steps'' and ``stages'' of Nash--Kuiper \cite{nash, k1, k2}, whose limit as $q \nearrow \infty$ yields the desired irregular isometric embedding. The following estimates are given in \cite[\S 1.6.1]{cds}:
\begin{eqnarray}
&&\max_{i,j}\|\p_i u_q \cdot \p_j u_q -g_{ij} \|_0 \leq \delta^2_q, \label{cds1} \\
&&\|u_q\|_2 \leq \mu_q,\label{cds2}\\
&&\left\|u_{q+1} - u_q\right\|_1 \leq C\delta_q,\label{cds3}
\end{eqnarray} 
where $\{\delta_q\}$ decreases exponentially while $\{\mu_q\}$ increases exponentially, provided that $\delta_0$ has been  appropriately chosen.  \eqref{cds2} implies that any limit of $u_q$ fails to be $C^2$-regular. To get a $C^{1,\alpha}$-limit with the range of $\alpha$ as stated in Theorem~\ref{thm: cds}, the authors of \cite{cds} made a delicate choice of $\delta_q, \mu_q$; for our purpose below, we only need the existence of such  parameters (but not the explicit choices) as in \eqref{cds1}--\eqref{cds3}.

We propose the following renormalisation process.

Let $\{u_q\}$ be the sequence of  smooth embeddings constructed  in Theorem \ref{thm: cds}. They are not isometric to the prescribed metric $g$ on $M$ in general (in fact, they are ``short''), and they take values in $\R^3$. We shall reserve the symbol $\mathfrak{e}$ for the Euclidean metric. Our proposed renormalisation process consists of the following ingrediences:
\begin{enumerate}
\item
Set $\gq := u_q ^\# \mathfrak{e}$, the pullback metric; assume that $u_q$ are $1$-periodic.
\item
Let $\kq$ be the Gauss curvature of $\gq$, and let $\kqone, \kqtwo$ be the principal curvatures.
\item
Let $\Hq$ be the second fundamental form of $u_q: \left(M, \gq\right) \map \left(\R^3,\mathfrak{e}\right)$.
\item
Let $(M,g)$ be the flat $2$-torus, which serves as the tentative limiting manifold.
\end{enumerate}

By construction, $u_q$ are isometric embeddings $\left(M, \gq\right) \emb (\R^3,\mathfrak{e})$. Thus $\left(\Hq, \kq\right)$ satisfies the Gauss--Codazzi equations with respect to $\gq$. Introduce the following new variables:
\begin{equation}\label{renorm variables}
\hq:=\frac{\Hq}{\eta_q},\qquad \gaq := \frac{\kq}{\eta_q^2},\qquad \gamma^{(q)}_i := \frac{\kappa^{(q)}_i}{\eta_q} \quad \text{for } i \in\{1,2\},
\end{equation}
where as in \eqref{cds2}, 
\begin{equation}\label{renorm parameter}
\eta_q := \|u_q\|_2 \leq \mu_q.
\end{equation}

To proceed, dividing by $\eta_q^2$ on both sides of \eqref{gauss}, we obtain that 
\begin{equation}\label{renorm gauss}
h^{(q)}_{11}h^{(q)}_{22} - h^{(q)}_{12}h^{(q)}_{21} = \det \gq \gaq.
\end{equation}
From the Codazzi equation \eqref{codazzi} one deduces that
\begin{align*}
\p_i H_{jk} - \p_jH_{ik} =  \G^l_{ik}H_{jl} - \G^l_{jk} H_{il},
\end{align*}
utilising the identity $\na_k \phi_{ij} = \p_k\phi_{ij} - \G^l_{ik}\phi_{lj} - \G^l_{jk}\phi_{il}$ for covariant $2$-tensor $\phi=\{\phi_{ij}\}$. Thus 
\begin{equation}\label{renorm codazzi}
\p_i h^{(q)}_{jk} - \p_j h^{(q)}_{ik} = ^{(q)}\G^l_{ik}h^{(q)}_{jl} - ^{(q)}\G^l_{jk} h^{(q)}_{il} := Q_{i,j,k}^{(q)}.
\end{equation}
Here $^{(q)}\G^i_{jk}$ are the Christoffel symbols of the Levi-Civita connection of $\left(M, \gq\right)$, thus $$\left|Q_{i,j,k}^{(q)}\right| \lesssim \mathcal{O}(\eta_q).$$ The process of Nash--Kuiper iterations can be viewed as defined 1-periodically in the (global) co-ordinates $(x_1, x_2) \in \R^2$. This crucially relies on our choice that $(M,g)$ is the flat torus.

Let us introduce the following change of the variables. We set
\begin{equation*}
z^{(q)}_j := \eta_q x_j
\end{equation*}
and
\begin{equation*}
\widehat{h^{(q)}_{ij}} \left(z^{(q)}_1,z^{(q)}_2\right)  := h^{(q)}_{ij} (x_1, x_2)
\end{equation*}
for $i, j \in \{1,2\}$. As a consequence, we may further express the renormalised Gauss and Codazzi equations~\eqref{renorm gauss} and \eqref{renorm codazzi} as
\begin{eqnarray}
&& \widehat{h^{(q)}_{11}}\,\widehat{h^{(q)}_{22}} - \widehat{h^{(q)}_{12}}\,\widehat{h^{(q)}_{21}} = \det\gq\, \gaq,\label{renorm gauss'}\\
&& \frac{\p}{\p z^{(q)}_i} \widehat{h^{(q)}_{jk}} - \frac{\p}{\p z^{(q)}_j} \widehat{h^{(q)}_{ik}} = \frac{Q_{i,j,k}^{(q)}}{\eta_q}. \label{renorm codazzi'}
\end{eqnarray}
Here and hereafter we have assumed a slight abuse of notations: we designate 
\begin{eqnarray*}
&&\gq\left(z^{(q)}_1,z^{(q)}_2\right)  \equiv \gq (x_1, x_2),\\
&&\gaq\left(z^{(q)}_1,z^{(q)}_2\right)  \equiv \gaq (x_1, x_2),\\
&&Q_{i,j,k}^{(q)}\left(z^{(q)}_1,z^{(q)}_2\right)  \equiv Q_{i,j,k}^{(q)} (x_1, x_2),
\end{eqnarray*}
and the like. Notice that $\widehat{h^{(q)}}$ is defined on $\R^2$ with period $\eta_q$; so in \eqref{renorm codazzi'} we may drop the superscript $^{(q)}$ in the derivatives. That is, one has  $\p_i \widehat{h^{(q)}_{jk}} - \p_j \widehat{h^{(q)}_{ik}} = {Q_{i,j,k}^{(q)}}\slash {\eta_q}$.

Consider now the renormalised equations \eqref{renorm gauss'} and \eqref{renorm codazzi'}. The variables $\widehat{\hq}, \gaq, \gq$ are uniformly bounded in $q$, thanks to the estimates and definitions in  \eqref{cds1}, \eqref{renorm variables}, and \eqref{renorm parameter}. Applying the theory of {\em compensated compactness} to the Gauss--Codazzi system ({\it cf.} Chen--Slemrod--Wang \cite{csw} and Chen--Li \cite{cl, clnew}, based on  foundational works of Murat \cite{m1, m2} and Tartar \cite{t1, t2}; \emph{cf}. also Dafermos \cite{dafermos}), we may  deduce the existence of measures $\overline{h}_{ij}$ and  $\overline{\gamma}$ satisfying
\begin{eqnarray}
&&\widehat{h^{(q)}_{ij}} \wstar \overline{h}_{ij},\qquad \gaq \wstar \overline{\gamma},\qquad \frac{Q_{i,j,k}^{(q)}}{\eta_q} \wstar F_{i,j,k}\quad \text{as } q \nearrow \infty; \label{wstar}\\
&&\p_1\lh_{22}-\p_2\lh_{12}=F_1,\qquad -\p_1\lh_{12} + \p_2 \lh_{11} = F_2; \label{bar-codazzi}\\
&&\lh_{11}\,\lh_{22} - \lh_{12}\,\lh_{21} = \overline{\gamma}.\label{bar-gauss}
\end{eqnarray}
Here \eqref{wstar} is understood in the weak-$\star$ topology of $L^\infty (\R^2)$, after passing to subsequences if necessary. The source term $F_{i,j,k}$ is given by
\begin{equation*}
F_{i,j,k} = \begin{cases}
F_1\qquad \text{ if } (i,j,k) = (1,2,2),\\
F_2 \qquad \text{ if } (i,j,k) = (1,2,1),\\
0\qquad \text{ if else}. 
\end{cases}
\end{equation*}
The equalities \eqref{bar-codazzi} and \eqref{bar-gauss} are understood in the sense of distributions and {\em a.e.}; they follow from the normalised equations  \eqref{renorm codazzi} and \eqref{renorm gauss}, as well as the estimate \eqref{cds1} obtained in \cite{cds}. 

\begin{lemma}\label{lem: renormalisation}
In the above setting, we have $F_1 = F_2 = \overline{\gamma} = 0.$
\end{lemma}

\begin{proof}

By the renormalised Codazzi equations \eqref{renorm codazzi} and \eqref{renorm codazzi'} we have
\begin{equation}\label{add A}
\p_i \widehat{h^{(q)}_{jk}} - \p_j \widehat{h^{(q)}_{ik}} =  \frac{^{(q)}\G^l_{ik}}{\eta_q} \widehat{h^{(q)}_{jl}} - \frac{^{(q)} \G^l_{jk}}{\eta_q} \widehat{h^{(q)}_{il}} \equiv \frac{Q_{i,j,k}^{(q)}}{\eta_q}.
\end{equation}
 As $\left|^{(q)}\G^i_{jk}\right| \lesssim \mathcal{O}(\eta_q)$, this quantity is uniformly bounded in $L^\infty(\R^2)$. Thus, invoking the Sobolev embedding lemma and relabelling the indices, $\left\{\p_i \widehat{h^{(q)}_{jl}} - \p_j \widehat{h^{(q)}_{il}}\right\}$ lies in a compact subset of $H^{-1}(\R^2)$. On the other hand, applying the same arguments as for the renormalised Gauss--Codazzi equations~\eqref{renorm gauss'} , \eqref{renorm codazzi'} and recalling the definition of Riemann curvature, we may redefine the independent variables  (by periodicity) to deduce that
\begin{equation}\label{add B}
\p_j \left(\frac{^{(q)}\G^l_{ik} }{\eta_q}\right) - \p_i \left(\frac{^{(q)}\G^l_{jk}}{\eta_q}\right) = \frac{^{(q)} R^l_{kji} }{\eta_q^2} - \frac{^{(q)}\G^l_{jp}}{\eta_q} \frac{^{(q)}\G^p_{ik}}{\eta_q} + \frac{^{(q)}\G^l_{ip}}{\eta_q} \frac{^{(q)}\G^p_{jk}}{\eta_q}.
\end{equation}
Here $^{(q)} R^l_{kji}$ is the Riemann curvature tensor of the Riemannian manifold $\left(M, u^\#_q\mathfrak{e}\right)$. This quantity is also confined in a compact subset of $H^{-1}(\R^2)$.

Therefore, by a standard application of the div-curl lemma (see \cite{m1,m2,t1,t2}) using the differential constraints  \eqref{add A} and \eqref{add B}, one may infer that
\begin{equation}\label{add C}
F_{i,j,k} = \overline{\G^l_{ik}}\overline{h_{jl}} - \overline{\G^l_{jk}}\overline{h_{il}}, 
\end{equation}
where $\overline{\G^i_{jk}}$ is the $L^\infty$-weak-$\star$ limit of $^{(q)}\G^i_{jk}\slash\eta_q$. Moreover, by definition one has
\begin{equation}\label{add D}
^{(q)}\G^i_{jk} = \frac{1}{2}\,{g^{(q)}}^{ip} \left\{ \p_j g^{(q)}_{pk} + \p_k g^{(q)}_{pj} - \p_p g^{(q)}_{jk} \right\}.
\end{equation}
Then, taking an arbitrary test function $\Phi \in C^\infty_0(\R^2)$ and  integrating by parts with respect to the Lebesgue measure, we get
\begin{align*}
\frac{1}{\eta_q} \int_{\R^2} {\gq}^{ip} \left(\p_j g^{(q)}_{pk}\right) \Phi\,\dd x &= \frac{1}{\eta_q}\int_{\R^2} \Phi \delta^{ip} \left(\p_j g^{(q)}_{pk}\right)\,\dd x + \frac{1}{\eta_q} \int_{\R^2}\Phi \Big\{{g^{(q)}}^{ip} -\delta^{ip}\Big\} \left(\p_j g^{(q)}_{pk}\right)\,\dd x\\
&= -\frac{1}{\eta_q} \int_{\R^2} \left(\p_j \Phi\right) g^{(q)}_{ik} \,\dd x + \int_{\R^2} \Big\{{g^{(q)}}^{ip} -\delta^{ip}\Big\}  \frac{\Phi \left(\p_i g^{(q)}_{pk}\right)}{\eta_q}\,\dd x \\
&=: {\rm I} + {\rm II}.
\end{align*}
Here ${\rm I} \map 0$ since the integral $\int_{\R^2} \left(\p_j \Phi\right) g^{(q)}_{pk}\,\dd x$ is uniformly bounded while $\eta_q \nearrow \infty$. In addition, ${\rm II} \map 0$ since ${g^{(q)}}^{ip} -\delta^{ip} \to 0$ in the $C^0$-topology while $\left\{\eta^{-1}_q \Phi\left(\p_i g_{pk}^{(q)}\right)\right\}$ is uniformly bounded. Applying similar arguments to the two other terms on the right-hand side of \eqref{add D}, we  deduce that
\begin{equation}\label{add E}
\int_{\R^2} \frac{^{(q)}\G^i_{jk}}{\eta_q}\, \Phi\,\dd x \longrightarrow 0\qquad \text{as } q \nearrow \infty.
\end{equation}
It thus follows that $\overline{\G^{i}_{jk}} \equiv 0$, hence $F_{i,j,k}\equiv 0$ for all $i,j,k \in \{1,2\}$. 

To see $\overline{\gamma} = 0$,  
 notice that $\overline{\gamma}$ is equal to the distributional limit of 
\begin{equation*}
\gaq = \frac{1}{\eta_q^2} g^{(q)}_{1j} \left\{ \p_1    {}^{(q)}\G^j_{22} - \p_2  {}^{(q)}\G^j_{12} + {}^{(q)}\G^j_{1k}  {}^{(q)}\G^k_{22} -   {}^{(q)} \G^j_{2k} {}^{(q)}\G^k_{12} \right\}.
\end{equation*}
The first term on the right-hand side converges to zero as $q\nearrow \infty$ in the sense of distributions, since for any test function $\Phi \in C^\infty_0(\R^2)$, it holds that
\begin{align*}
\frac{1}{\eta_q^2} \int_{\R^2}\Phi g^{(q)}_{1j} \left(\p_1 {}^{(q)}\G^j_{22}\right)\,\dd x &=-\frac{1}{\eta_q^2} \int_{\R^2} \p_1\Phi g^{(q)}_{1j} \G^j_{22}\,\dd x - \frac{1}{\eta_q^2} \int_{\R^2} \Phi \p_1 g^{(q)}_{1j} \G^j_{22} \,\dd x\\
& \longrightarrow 0,
\end{align*}
by the definition of $\eta_q$ and the identity $\overline{\G}^i_{jk} \equiv 0$ established above. The same argument applies to the second term. As a consequence, we have that
\begin{equation}\label{gamma-bar}
\overline{\gamma} = \lim_{q\map\infty} \frac{{}^{(q)}\G^1_{1k}  {}^{(q)}\G^k_{22} -   {}^{(q)} \G^1_{2k} {}^{(q)}\G^k_{12}}{\eta^2_q},
\end{equation} 
where the limit is understood in the sense of distributions. 

Now we apply once again the classical div-curl lemma (see  \cite{m1,m2,t1,t2}). Consider the vector fields on $\R^2$:
\begin{equation*}
V_q:=\eta_q^{-1}\left({}^{(q)}\G^1_{1k},\,{}^{(q)}\G^1_{2k}\right)^\top
\end{equation*}
and 
\begin{equation*}
W_q:=\eta_q^{-1}\left({}^{(q)}\G^k_{22}, \,-{}^{(q)}\G^k_{12}\right)^\top.
\end{equation*}
The curl of $V_q$ and the divergence of $W_q$ are uniformly bounded in $L^\infty(\R^2)$. Thus, as $\overline{\G^i_{jk}}\equiv 0$, we have
\begin{equation*}
\overline{\gamma} = \overline{\G^1_{1k}} \, \overline{\G^k_{22}} - \overline{\G^1_{2k}} \,\overline{\G^k_{12}} = 0.
\end{equation*}

Hence the lemma is proved.  \end{proof}


By virtue of \eqref{bar-codazzi}, \eqref{bar-gauss} and Lemma~\ref{lem: renormalisation}, we define as in \S\ref{sec: 2d} the fluid variables $\left\{\rho, v^1, v^2, p\right\}$ via the following relations:
\begin{equation}\label{renorm: rho, p v by h-bar}
\begin{cases}
\rho \left(v^1\right)^2 + p = \lh_{22},\\
\rho v^1 v^2 = -\lh_{12},\\
\rho \left(v^2\right)^2 + p = \lh_{11}.
\end{cases}
\end{equation}
The limiting equations \eqref{renorm codazzi'} and \eqref{renorm gauss'} thus become
\begin{equation}\label{renorm + gc}
\begin{cases}
\p_1 \left(\rho \left(v^1\right)^2 + p\right) + \p_2 \left(\rho v^1 v^2\right) = 0,\\
\p_1 \left(\rho v^1v^2\right) + \p_2\left(\rho \left(v^2\right)^2 + p\right)= 0,\\
\rho p\left[ \left(v^1\right)^2+\left(v^2\right)^2\right] + p^2 = 0.
\end{cases}
\end{equation}
With $\lh$ given as above, we may solve for $p$ in terms of $\lh$. Indeed, the first and the third equations of \eqref{renorm: rho, p v by h-bar} yield that $$p\left[ \left(v^1\right)^2+\left(v^2\right)^2\right]  +2p = \lh_{11}+\lh_{22}.$$ Together with \eqref{renorm gauss} and the last equation in \eqref{renorm + gc}, one deduces that
\begin{align}\label{p: sol, renorm}
p = 0 \qquad \text{ or }\qquad p = \overline{h}_{11} + \overline{h}_{22}.
\end{align}
With either choice of $p$, as in \S \ref{sec: 2d} we may solve for $v^1$, $v^2$, and $\rho$ in the following way. Denote
\begin{equation*}
\overline{f_{11}} := \lh_{11} - p,\qquad \overline{f_{22}} := \lh_{22} - p.
\end{equation*}
Since $\overline{f_{11}}\,\overline{f_{22}} - \lh_{12}\lh_{21}=0$ (where $\lh$ is symmetric), without loss of generality we may assume $\overline{f_{11}} >0$ and $\overline{f_{22}} >0$. Then,  we require $\rho >0$ to satisfy the continuity equation
\begin{equation}\label{renorm: continuity equation}
\p_1 \sqrt{\rho \overline{f_{11}}} + \p_2 \sqrt{\rho \overline{f_{22}}}  = 0,
\end{equation}
and solve $v$ from
\begin{equation}\label{renorm: v eq}
v_1 = \left(\frac{\overline{f_{11}}}{{\rho}}\right)^{\frac{1}{2}},\qquad v_2 =\left( \frac{\overline{f_{22}}}{{\rho}}\right)^{\frac{1}{2}}.
\end{equation}

In terms of $\overline{h}_{ij}$, \eqref{renorm: continuity equation} and \eqref{renorm: v eq} are equivalent to
\begin{equation}\label{x1}
\p_1 \sqrt{\rho \overline{h}_{11}} + \p_2 \sqrt{\rho \overline{h}_{22}} = 0\qquad \text{ or } \qquad \p_1\sqrt{-\rho h_{22}} + \p_2 \sqrt{-\rho h_{11}}=0
\end{equation}
and
\begin{equation}\label{x2}
(v_1,v_2) =\left( \left[\frac{\overline{h}_{11}}{\rho}\right]^{\frac{1}{2}}, \left[\frac{\overline{h}_{22}}{\rho}\right]^{\frac{1}{2}}\right)\quad \text{or}\quad (v_1,v_2) =\left( \left[-\frac{\overline{h}_{22}}{\rho}\right]^{\frac{1}{2}}, \left[-\frac{\overline{h}_{11}}{\rho}\right]^{\frac{1}{2}}\right)
\end{equation}
with respect to the choices of the two roots of $p$ in \eqref{p: sol, renorm}.

If $h_{11}, h_{22} >0$, we choose the pressureless case $p=0$, and if $h_{11}, h_{22} < 0$, we choose $p= \overline{h}_{11} + \overline{h}_{22}$. Notice that the two solutions in \eqref{x1} and \eqref{x2} can be transformed into each other by identifying $$\left(\overline{h}_{11}, \overline{h}_{22}\right)\longmapsto \left(-\overline{h}_{22}, -\overline{h}_{11}\right).$$ We may thus regard them as dynamically indistinguishable. That is, both choices in \eqref{p: sol, renorm} describe the motion of a pressureless gas. Furthermore, if $\overline{h}_{11}, \overline{h}_{22}$ are smooth, we may use the method of characteristics to  locally define $\rho$ and thus accomplish our
goal of delivering the fluid variables $\{\rho, v, p\}$. In general, however, we only have $\overline{h}_{11}, \overline{h}_{22}\in L^\infty(M)$, so the existence (even definition) of the solution is not obvious.


To summarise this section, we have derived the following
\begin{theorem}\label{thm: renormalization}
The formulae \eqref{p: sol, renorm}, \eqref{renorm: continuity equation}, and \eqref{renorm: v eq} provide expressions for the fluid variables $\{p, v,\rho\}$ of the renormalised limit of the Nash--Kuiper iterations. They describe a weak solution to the steady  pressureless gas equations, under the assumption that a solution $\rho$ exists to the continuity equation \eqref{x1}. 
\end{theorem}

\section{Extension to Multi-Dimensions}\label{sec: multi D}

Now we aim at extending our earlier investigations to multi-dimensions. That is, we hope to establish a link between the isometrically embedded smooth submanifolds in $\R^{n+k}$ (for arbitrary $n, k$) and the solutions to steady compressible Euler equations.

Let $(\Sigma,g)$ be an $n$-dimensional Riemannian manifold isometrically embedded in $\R^{n+k}$. The extrinsic geometry is characterised by the second fundamental from $\left\{H^\mu_{ij}\right\}$ and the affine  connection on the normal bundle $\left\{A^\nu_{\mu i}\right\}$. In this section, $i,j,k,l \in \{1,\ldots,n\}$ are  indices for the tangent bundle $T\Sigma$, and $\nu,\mu,\eta \in\{n+1,\ldots,n+k\}$ are indices for the normal bundle $T^\perp\Sigma$. In this case, the fundamental equations for the existence of isometric embeddings are the Gauss, Codazzi and Ricci equations as follows:
\begin{eqnarray}
&& H^\mu_{ij}H^\mu_{kl} - H^\mu_{ik} H^\mu_{jl} = R_{iljk},\label{gauss nd}\\ 
&& \na_i H^{\mu l}_j - \na_j H^{\mu l}_i + A^\nu_{\mu i}H^{\nu l}_j - A^\nu_{\mu j} H^{\nu l}_i = 0,
\label{codazzi nd}\\
&& \na_i A^\nu_{\mu j} - \na_j A^\nu_{\mu i} + A^\nu_{\eta i} A^\eta_{\mu j} - A^{\nu}_{\eta j}A^\eta_{\mu i} = g^{pq}\big(H^\mu_{ip}H^\nu_{jq}-H^\mu_{jp}H^\nu_{iq}\big).\label{ricci nd}
\end{eqnarray}
The mean curvature vector $$\vec{\mean}=\left(\mean^1, \ldots, \mean^k\right)$$ is given by
\begin{equation*}
\mean^\mu := g^{ij}H^\mu_{ij} =H^{\mu i}_i. 
\end{equation*}

On $(\Sigma,g)$ we also have the steady compressible Euler equation. The continuity equation is
\begin{equation}\label{continuity eq nd}
\na_k \left(\rho v^k\right) = 0.
\end{equation}
Moreover, for the stress-energy tensor
\begin{equation}\label{def of stress energy tensor, nd}
P^{ij} = \rho v^i v^j + g^{ij}p,
\end{equation}
there holds the balance of linear momentum:
\begin{equation}\label{P eq nd}
\na_j P^j_i = \Pi_i,
\end{equation}
The unknowns of the Euler equations are the fluid variables $\rho$, $v=\left\{v^i\right\}_{i=1}^n$, and $p$. Here $\Pi_i$ is a given body force on $\Sigma$. 

Assume that $H_{ij}$ is a smooth solution for the Gauss--Codazzi--Ricci equations \eqref{gauss nd}--\eqref{ricci nd}. We shall  identify the fluid variables as functions of $H_{ij}$. 

To begin with, consider the {\em contracted Codazzi equation}:
\begin{equation*}
\na_i H^{n+1, j}_j - \na_j H^{n+1, j}_i = \Pi_i,
\end{equation*}
where
\begin{equation*}
\Pi_i  := A^\nu_{n+1, j}H^{\nu j}_i - A^{\nu}_{n+1, i}H^{\nu j}_j.
\end{equation*}
This is obtained by setting $\mu = n+1$ and multiplying by $\delta^j_l$ in \eqref{codazzi nd}. We  choose the following relation between the stress-energy tensor and second fundamental forms:
\begin{equation}\label{P def, nD}
P^{j}_i := -H^{n+1, j}_i + \mean^{n+1} \delta^j_i.
\end{equation}
Then
\begin{align*}
\na_j P^j_i = -\na_i H^{n+1, j}_j + \na_i \mean^{n+1} +\Pi_i = \Pi_i.
\end{align*}
{\it i.e.}, the balance law for linear momentum \eqref{P eq nd} is satisfied.


By contracting with the metric tensor, we see that \eqref{P def, nD} is equivalent to 
\begin{equation}
H^{n+1}_{ij} = \mean^{n+1} g_{ij} - P_{ij}. 
\end{equation}
We can deduce the constitutive relation for $P_{ij}$ from the Gauss equation \eqref{gauss nd}. Indeed, define
\begin{equation}
L_{iljk} := \sum_{\mu =n+2}^{n+k}  H^\mu_{ij}H^\mu_{kl} - H^\mu_{ik} H^\mu_{jl}
\end{equation}
Then one may compute
\begin{align*}
R_{iljk} &= (-P_{ij} +\mean^{n+1} g_{ij})(-P_{kl} + \mean^{n+1} g_{kl}) \\
&\qquad\qquad- (-P_{ik}+\mean^{n+1} g_{ik}) (-P_{jl}+\mean^{n+1} g_{jl})+ L_{iljk}\\
&= \Big\{P_{ij}P_{kl} - P_{ik} P_{jl}\Big\} -\mean^{n+1} \Big\{ P_{kl}g_{ij} + P_{ij}g_{kl} - P_{ik}g_{jl} - P_{jl} g_{ik} \Big\}\\
&\qquad\qquad  + [\mean^{n+1}]^2\Big\{ g_{ij}g_{kl} - g_{ik}g_{jl} \Big\}+ L_{iljk}.
\end{align*}
Substituting \eqref{def of stress energy tensor, nd} into the above expression, one obtains that
\begin{align*}
R_{iljk} &= \Big\{(\rho v_iv_j+pg_{ij})(\rho v_kv_l + pg_{kl}) - (\rho v_iv_k + pg_{ik})(\rho v_jv_l + pg_{jl})\Big\}\\
&\quad -\mean^{n+1} \Big\{ g_{ij} (\rho v_kv_l + pg_{kl}) + g_{kl}(\rho v_iv_j + pg_{ij})  - g_{jl} (\rho v_iv_k + pg_{ik}) - g_{ik} (\rho v_jv_l + pg_{jl}) \Big\}\\
&\quad + [\mean^{n+1}]^2 \Big\{ g_{ij}g_{kl}-g_{ik}g_{jl} \Big\}+ L_{iljk}\\
&= (p-\mean) A_{iljk} + (p-\mean)^2 B_{iljk} + L_{iljk},
\end{align*}
where
\begin{equation}\label{A,B def in nd}
\begin{cases}
A_{iljk} = \rho\Big[v_iv_j g_{kl} + v_kv_lg_{ij} - v_iv_k g_{jl}-v_jv_lg_{ik}\Big],\\
B_{iljk} = g_{ij}g_{kl} - g_{ik}g_{jl}. 
\end{cases}
\end{equation}
These terms have a natural geometric structure: for $2$-tensors $T=T_{ij}$ and $S=S_{ij}$, denote by $T\odot S$ the $4$-tensor $$(T\odot S)_{iljk}:=T_{ij}S_{kl}-T_{ik}S_{jl}.$$ Then
\begin{equation*}
\begin{cases}
A = 2\rho(\sigma \odot \sigma),\\
B=g\odot g,
\end{cases}
\end{equation*}
where $$\sigma := \frac{1}{2}(g\otimes (v\otimes v) + (v\otimes v) \otimes g)$$ is the symmetrisation of $g$ and $v\otimes v$.

Now we express $p$ in terms of the geometric quantities. This is done by contracting the Riemann curvature. First, we compute the Ricci curvature tensor $\ric$:
\begin{align}
\ric_{lk} &= g^{ij}R_{iljk} \nonumber\\
&= (p-\mean) g^{ij} A_{iljk} + (p-\mean)^2 g^{ij} B_{iljk} +  g^{ij}L_{iljk}\nonumber\\
&= (p-\mean) \rho \Big[ (n-2)v_kv_l + g_{kl}|v|_g^2\Big] +(p-\mean)^2(n-1)g_{kl} +  g^{ij}L_{iljk}.
\end{align}
Here $|v|_g^2 := g^{ij} v_iv_j$. Contracting once more yields
\begin{align}
\scal &= g^{kl}\ric_{kl}\nonumber\\
&= (n-1)(p-\mean) \Big[ 2\rho|v|^2_g + n(p-\mean)\Big] +s,
\end{align}
where
\begin{equation}
s:=g^{ij}g^{kl}L_{iljk}.
\end{equation}
Thanks to \eqref{def of stress energy tensor, nd} and \eqref{P def, nD}, the momentum energy density can be expressed in terms on $p$ and $\mean$:
\begin{align*}
\rho|v|^2_g &= P^{n+1,i}_i - np \\
&= -H^{n+1,i}_i + n\mean^{n+1} - np \\
&= (n-1)\mean^{n+1} - np.
\end{align*}

Now we can conclude that $p$ satisfies the quadratic equation:
\begin{equation}\label{quadratic eq for p, nd}
-n(n-1)p^2 + 2(n-1)^2\mean^{n+1}p + (n-1)(n-2)[\mean^{n+1}]^2 - \scal+s =0.
\end{equation}
It has real roots
\begin{equation}\label{p roots nd}
p = \frac{(n-1)^2\mean^{n+1} \mp \sqrt{ (n-1)^4[\mean^{n+1}]^2 + n(n-1)^2(n-2) [\mean^{n+1}]^2 + n(n-1)(s-\scal) } }{n(n-1)},
\end{equation}
whenever the discriminant $$\Delta =  (n-1)^4[\mean^{n+1}]^2 + n(n-1)^2(n-2) [\mean^{n+1}]^2 + n(n-1)(s-\scal) \geq 0.$$ 

Similar to $\S \ref{sec: 2d}$, we can solve for the velocity from 
\begin{equation}
v^j = \sqrt{\frac{f^{jj}}{\rho}}
\end{equation}
(no summation convention), where  $\rho$ is determined by the continuity equation \eqref{continuity eq nd}, and
\begin{equation}
f^{ij} := -H^{n+1, i,j} + (\mean^{n+1}-p)g^{ij}. 
\end{equation}
This can be done whenever for all $q,k$ there holds
\begin{equation}\label{consistency}
\rho v_qv_k = -H^{n+1}_{qk} + (\mean^{n+1}-p)g_{qk}.
\end{equation}
By considering $\rho (v_q)^2 \rho (v_k)^2 = (\rho v_qv_k)^2$, we can re-express \eqref{consistency} as the following {\em consistency conditions} (C1), (C2): 
\begin{align*}
\text{ For all $k,q$, the ``principal matrices'' $G^{-1}C$ have a common eigenvalue $\lambda$},\tag{C1}
\end{align*}
where \begin{equation*}
G = \begin{bmatrix}
g_{kk} & g_{kq}\\
g_{qk} & g_{qq}
\end{bmatrix},\qquad C=\begin{bmatrix}
H^{n+1}_{kk} & H^{n+1}_{kq}\\
H^{n+1}_{qk} & H^{n+1}_{qq}
\end{bmatrix},
\end{equation*}
as well as
\begin{equation*}
\lambda = \mean^{n+1}-p \text{ with $p$ satisfying 
Equation \eqref{p roots nd}}. \tag{C2} 
\end{equation*}

At this point we could repeat the renormalisation procedure in \S \ref{sec: renorm}, when $M$ is the flat $n$-torus. The arguments would be similar to those in \S \ref{sec: renorm}, but now applied to the Gauss--Codazzi--Ricci equations in higher dimensions and codimensions ({\it cf. e.g.}, Chen--Slemrod--Wang \cite{csw2} and Chen--Li \cite{cl}): this is because the quadratic terms $H \otimes H$, $A \otimes A$, $H \otimes A$, $\G \otimes A$, and $\G \otimes H$ are of the order $\mathcal{O} (\eta_q^2)$, where $\eta_q$ is the $C^2$-norm of the $q$-th approximate solution $u_q$, and $\G$ is the Levi-Civita connection of $\left(M, u_q^\# \mathfrak{e}\right)$ as before. However, the consistency conditions $(C1)\,\&\,(C2)$ would still have to be satisfied for the weak-$\star$ limit equations.

As the closing remark,  correspondences between solutions to fluid dynamical PDE and solutions to geometric PDE (assuming sufficient regularity) have been studied in recent physics literature. When the geometric PDE are the vacuum Einstein equations for  a Ricci-flat hypersurface, by Damour \cite{damour} and Bredberg--Strominger \cite{phys4}, rather amazingly, the fluid dynamic PDE formally are the classical non-steady Navier--Stokes equations of incompressible fluids. See also \cite{phys1, phys2, phys3, phys5, phys6, phys7} and the references cited therein; this list is by no means exhaustive. It is interesting to further explore such geometry-fluid correspondences using analytic methods.

\begin{appendices}

\section{Geometry of the standard torus $\tor$}\label{sec: app A}

In this appendix, we collect some basic facts about the standard torus, which are used in \S\ref{subsec: torus}. The standard torus $\tor$ of major radius $a$ and minor radius $b$ is the parametric surface in $\R^3$ given by \eqref{torus}, reproduced below:
\begin{equation*}
\tor = \left\{\begin{bmatrix}
(a+b\cos \theta)\cos \phi\\
(a+b\cos \theta)\sin \phi\\
b\sin\theta
\end{bmatrix}:\,0 \leq \theta,\phi < 2\pi\right\}.
\end{equation*}
More generally, we  also consider the scaled version:
\begin{align}\label{torc, appendix}
\torc  = \left\{\begin{bmatrix}
\left[a+b\cos \left(\frac{\theta}{c}\right)\right]\cos \phi\\
\left[a+b\cos \left(\frac{\theta}{c}\right)\right]\sin \phi\\
b\sin\left(\frac{\theta}{c}\right)
\end{bmatrix}:\,0 \leq \phi < 2\pi,\,0 \leq \theta < 2\pi c\right\}\quad \text{where } c>0.
\end{align}
Note that $\mathbf{T}_1(a,b)=\tor$.

The Riemannian metric (\emph{i.e.}, the first fundamental form) of $\torc$ is
\begin{equation}\label{metric, tor}
g = \begin{bmatrix}
\left[a+b \cos \left( \frac{\theta}{c} \right)\right]^2 & 0\\
0&\frac{b^2}{c^2}
\end{bmatrix},
\end{equation}
and the second fundamental form $H=\{H_{ij}\}$ is
\begin{equation}\label{second, tor}
H = \begin{bmatrix}
-\left[a+b \cos \left( \frac{\theta}{c} \right)\right]\cos \left( \frac{\theta}{c} \right)  & 0\\
0& -\frac{b}{c^2}
\end{bmatrix}.
\end{equation}
From these we can compute the Gauss curvature
\begin{align}\label{gauss curv, tor}
\kappa &= \frac{\det H}{\det g} = \frac{\cos \left( \frac{\theta}{c} \right)}{b\left[a+b \cos \left( \frac{\theta}{c} \right)\right]}
\end{align}
and the mean curvature 
\begin{align}\label{mean curv, tor}
\mean &= \frac{H_{11}g_{22} - 2H_{12}g_{12} + H_{22}g_{11}}{\det g} \nonumber \\
&= - \frac{a+2b \cos \left( \frac{\theta}{c} \right) }{b\left[a+b \cos \left( \frac{\theta}{c} \right)\right]}.
\end{align}
Then, the principal curvatures $\kappa_\pm$ are the roots for the quadratic polynomial:
\begin{align*}
q(s) := s^2 - \mean s + \kappa,
\end{align*}
and hence
\begin{align*}
s = -\frac{\cos  \left( \frac{\theta}{c} \right)}{\left[a+b \cos \left( \frac{\theta}{c} \right)\right]}\quad \text{ or } \quad -\frac{1}{b}.
\end{align*}

One issue arises here: we have chosen $p=\kappa_-$, the smaller principal curvature in \eqref{p=kappa 2, 2d}.  But the smaller solution to $q(s)=0$ above is $s=-1/b$, while it is non-physical to keep the pressure fixed for a compressible fluid. Moreover, throughout \S\ref{subsec: torus} we have assumed the $\kappa_+$ is constant.  This issue is resolved as in Shiohama--Takagi \cite[p.479]{st}: by inverting the orientation we replace the unit normal vector field $e_3$ by $-e_3$. In this way, \eqref{second, tor} and \eqref{mean curv, tor} are replaced by 
\begin{equation}\label{second, mean, corrected}
H = \begin{bmatrix}
\left[a+b \cos \left( \frac{\theta}{c} \right)\right]\cos \left( \frac{\theta}{c} \right)  & 0\\
0& \frac{b}{c^2}
\end{bmatrix}\quad \text{ and } \quad \mean = \frac{a+2b \cos \left( \frac{\theta}{c} \right) }{b\left[a+b \cos \left( \frac{\theta}{c} \right)\right]}.
\end{equation} 
The principal curvatures are 
\begin{equation}\label{ppl curv, tor}
\kappa_+ = \frac{1}{b}\quad \text{ and } \quad \kappa_- = \frac{\cos  \left( \frac{\theta}{c} \right)}{\left[a+b \cos \left( \frac{\theta}{c} \right)\right]} = p.
\end{equation}
Also note that our convention for mean curvature is $\mean = {\rm tr}_g H$ instead of $\mean = \frac{1}{\dim M} {\rm tr}_g H$. Thus, the round 2-sphere has $\mean =2$ instead of $1$.

As a result, the flow speed, which is equal to the difference between the principal curvatures modulo a constant, is
\begin{equation}\label{q, torc}
q = |v|_g = \kappa_+ - \kappa_- = \frac{a}{b\left[a+b \cos \left( \frac{\theta}{c} \right)\right]}.
\end{equation} 
The fluid density as chosen in \S\ref{subsec: torus} (also modulo a constant) is
\begin{align*}
\rho = \frac{1}{\kappa_+ - \kappa_- } = \frac{b\left[a+b \cos \left( \frac{\theta}{c} \right)\right]}{a}.
\end{align*}
See \eqref{torus sonic} and \eqref{rho solution}. For the purpose of solving the system~\eqref{system} for irrotationality, conservation of mass, and Bernoulli's law, we may choose as above   the normalisation constants for both $q$ and $\rho$ to be 1, thanks to the scaling of these PDE.

Now let us write the equation
\begin{align*}
{\rm div}_g \left(\frac{v}{|v|_g}\right) = 0
\end{align*}
in local co-ordinates  on $\torc$. The natural co-ordinate system $\{\p_\theta, \p_\phi\}$ comes from the parametrisation of $\torc$ in \eqref{torc, appendix}. Then for $v \in \G(T\torc)$ we may write
\begin{align}
v = v^\theta (\theta,\phi) \p_\theta + v^\phi (\theta,\phi) \p_\phi. 
\end{align}
Expressing the Riemannian divergence as
\begin{align*}
{\rm div}_g \left(\frac{v}{q}\right) = \frac{1}{\sqrt{\det g}} \p_i \left(\sqrt{\det g}\, \frac{v^i}{q} \right)
\end{align*}
with Einstein's summation convention over $i \in \{\theta,\phi\}$, we obtain
\begin{equation}\label{Eq1, 1-laplace in coordinates}
-\frac{2b}{c}\sin\left(\frac{\theta}{c}\right) v^\theta(\theta,\phi) + \left[ a+ b\cos \left(\frac{\theta}{c}\right) \right] \left\{\p_\theta v^\theta(\theta,\phi) + \p_\phi v^\phi(\theta,\phi)\right\} = 0.
\end{equation}
In addition, we have the constraint $|v|_g = q$, where $q$ is determined by the geometric parameters of the torus $\torc$ as in \eqref{q, torc}. This is expressed in local co-ordinates by
\begin{equation}\label{Eq2, v in coordinates}
\left[ a+ b\cos \left(\frac{\theta}{c}\right) \right]^2 \left(v^\theta(\theta,\phi)\right)^2 + \frac{b^2}{c^2}\left(v^\phi(\theta,\phi)\right)^2 = \frac{a^2}{b^2\left[ a+ b\cos \left(\frac{\theta}{c}\right) \right]^2}.
\end{equation}

\section{Irrotational Chaplygin gas on surfaces} \label{sec: appendix, variation}

In this appendix, we prove Theorem~\ref{thm: chaplygin} (reproduced below) concerning the Chaplygin gas on an arbitrary non-simply-connected closed Riemannian manifold $(M,g)$. As for aerodynamists' convention, the system of irrotationality, conservation of mass, and Bernoulli's law will be considered. 

\begin{theorem}
Let $(M,g)$ be any closed Riemannian manifold. Let $h \in \Omega^1(M)$ be any nontrivial harmonic 1-form. There exist $\rho: M \to [0,\infty[$ and $v \in L^1(M;TM)$ which are weak solutions to the steady Euler equations away from the set of stagnation points of $v$: 
\begin{equation}\label{chaplygin, appendix}
\begin{cases}
\omega = 0,\\
{\rm div}_g(\rho v) = 0,\\
\rho q = {\rm constant},
\end{cases}
\end{equation}
such that $v^\sharp$, the 1-form canonically dual to $v$, lies in the same cohomology class as $h$. 

If, in addition, $(M,g)$ admits an everywhere non-vanishing harmonic 1-form, then $v$ can be chosen without stagnation points.
\end{theorem}

Here, in arbitrary dimensions, the vorticity $\omega$ is defined as the 2-form:
\begin{equation*}
\omega = \dd \left(v^\sharp\right).
\end{equation*}
By writing $v \in L^1(M;TM)$ we mean that $v$ is a vector field on $M$ with $L^1$-regularity. In general, for a vector bundle $E \to M$, we write $X(M;E)$ for $X = L^p, W^{k,p}, C^{k,\alpha},\ldots$ to denote the space of $E$-sections with indicated regularity $X$. 

\begin{remark}\label{rem: stagnation}
 The case $M = \{q=0\}$ is ruled out by restricting to the nontrivial cohomology class $[h]\neq [0]$. In 2 dimensions such $h$ exists if and only if $M$ is not a topological 2-sphere, by Hodge decomposition and de Rham's theorem. See, \emph{e.g.}, Petersen \cite[\S 7.2, Theorem~47 on p.205, and Appendix, Theorem~89 on p.386]{p}.
\end{remark}


Note that if $(M,g)$ is the Euclidean space, then the irrotationality of $v$ implies that $v=\na\psi$ for a stream function $\psi$. Then \eqref{1-laplace} becomes the \emph{1-Laplace equation} for the scalar field $\psi$:
\begin{align*}
{\rm div}\left(\frac{\na \psi}{|\na \psi|}\right) = 0,
\end{align*} 
which is the Euler--Lagrange equation for the energy functional $I[\psi]:=\int |\na\psi|\,\dd x$. When $M$ is not simply-connected, such stream function does not exist in general. We shall make use of a variational argument  adapted from Evans \cite[\S 3.1]{evans}.


\begin{proof}
We divide our arguments into three steps.

\smallskip
\noindent
{\bf 1.} Our strategy is as follows. First, by the scaling invariance property of the system~\eqref{chaplygin, appendix} (see \eqref{scaling}), we may assume that the constant therein is 1.  Substituting the Bernoulli's law $\rho q= 1$ into the second equation, we see that \eqref{chaplygin, appendix} is equivalent to
\begin{equation}\label{1-laplace, appendix}
\begin{cases}
\omega = 0,\\
{\rm div}_g\left( \frac{v}{q} \right) = 0.
\end{cases}
\end{equation}


Denote by $\alpha := v^\sharp \in \Omega^1(M,g)$, the differential 1-form canonically isomorphic to $v$. That $v$ is irrotational is equivalent to $\dd \alpha =0$, namely that $\alpha$ is a  closed 1-form. Thus, for the given harmonic 1-form $h$, it is natural to solve for ${\rm div}_g\left( \frac{v}{q} \right) = 0$ in the cohomology class 
\begin{equation*}
[h] := \left\{ \beta \in \Omega^1(M,g):\, \beta - h = \dd\chi\text{ for some } \chi \in \Omega^0(M,g)\right\}. 
\end{equation*}
For regularity considerations, we shall relax to a larger class:
\begin{equation}\label{L1 class}
[[h]] := \left\{ \beta \in L^\infty(M;TM):\,\beta - h = \dd\chi \text{ in the sense of distributions for } \chi \in \mathcal{D}'(M)  \right\}.
\end{equation}
When $\beta \in [[h]]$ let us still say that $\beta$ is cohomologous to $h$ (in the setting of little regularity), as in Ciarlet--Gratie--Mardare \cite[\S 3]{cgm}.

For this purpose, the Hodge decomposition implies that there exist a scalar field $\psi$ and a harmonic 1-form $h \in \Omega^1(M,g)$  such that $\alpha = \dd \psi + h$. By raising and lowering indices, \eqref{1-laplace} is equivalent to
\begin{align}\label{EL}
\dd^* \left( \frac{\alpha}{\left|\alpha\right|_g} \right) = \dd^* \left( \frac{\dd \psi + h}{\left|\dd \psi + h\right|_g} \right) = 0.
\end{align}
It then remains to establish the following

\noindent
\emph{Claim A}: Fix any harmonic 1-form $h \in \Omega^1(M,g)$. There exists $\alpha \in L^\infty(M;T^*M)$ cohomologous to $h$ which is a weak solution to \eqref{EL}. The notion of weak solution is understood with its domain away from the zeros of $\alpha$ as in Remark~\ref{rem: stagnation}. That is, $\alpha \in [[h]]$ is a weak solution if and only if
\begin{align*}
\int_{M \setminus \{\alpha = 0\}} \bra \frac{\alpha}{|\alpha|_g},\,\dd\psi \ket_g \,\dvg = 0 \qquad \text{for any $\psi \in C^1_0\left(M \setminus \{\alpha = 0\}\right)$}.
\end{align*}

Indeed, assuming \emph{Claim A}, we may take $v = \alpha^\flat$ to solve for \eqref{1-laplace, appendix}, and then recover the Bernoulli law by setting
\begin{equation*}
\rho = \begin{cases}
\left(|v|_g\right)^{-1}\qquad \text{ when } v \neq 0;\\
\infty \qquad \text{ elsewhere}.
\end{cases}
\end{equation*}

As a side remark, for many surfaces, \emph{e.g.}, $M = \torc$ as in Appendix~\ref{sec: app A}, there exist everywhere non-vanishing harmonic 1-forms. Indeed, in the parametrisation \eqref{torc, appendix} for the torus, the space of harmonic vector fields is spanned by $\{\p_\theta,\,\p_\phi\}$. A linear combination $h= a \,\dd\theta + b\,\dd\phi$ is non-vanishing for generic $a,b\in\R$. Also, if $(M,g)$ has positive Gauss curvature, then any nontrivial harmonic 1-form has no zero by Bochner's theorem. See, \emph{e.g.}, Petersen \cite[\S 7.3, Theorem~48 and Corollary~18 on p.208]{p}.

\smallskip
\noindent
{\bf 2.} To prove \emph{Claim A} above, we resort to a variational formulation of \eqref{EL}, motivated by Evans \cite[\S 3.1]{evans}. Consider a subclass of $[[h]]$ defined in \eqref{L1 class}:
\begin{equation}\label{Ah}
\mathcal{A}_h := \left\{ \beta \in [[h]]:\, \|\beta\|_{L^\infty(M,g)} \leq 10 \|h\|_{L^\infty(M,g)} \right\}.
\end{equation}
We shall first establish

\noindent
\emph{Claim B}: \eqref{EL} is the Euler--Lagrange equations for the minimisation problem $\inf_{\mathcal{A}_h}\,  \mathcal{I}$, where
\begin{align}
\mathcal{I}[\beta]:= \|\beta\|_{L^1(M,g)}:=\int_M |\beta|_g \,\dvg.
\end{align}

\begin{proof}[Proof of Claim B] Let us compute the first variation of $\mathcal{I}$. Indeed, to move $\alpha$ within the cohomology class of $h$, we take perturbations of the form $\left\{\alpha + t \dd \eta\right\}_{-\delta_0 < t< \delta_0}$ for an arbitrary $\eta \in C^\infty(M)$ whose $C^1$-norm is bounded from the above by a uniform constant $c_0$, once $\delta_0$ is fixed. By Hodge decomposition, there is a scalar field $\psi$ such that $\alpha = \dd \psi + h$. Then
\begin{align*}
\frac{d}{dt} \mathcal{I}[\alpha + t\dd\eta] &= \int_M \left\{\frac{d}{dt} \left|h+\dd\psi + t\dd\eta\right|_g \right\}\,\dvg\\
&= \int_M \frac{2\bra h+\dd\psi,\,\dd\eta\ket_g + 2t|\dd\eta|^2_g}{\left|h+\dd\psi + t\dd\eta\right|_g } \,\dvg.
\end{align*}
So the minimality condition $\frac{d}{dt} \mathcal{I}[\alpha + t\dd\eta]\bigg|_{t=0} = 0$ leads to 
\begin{align*}
\int_M \bra\frac{ h+\dd\psi}{\,\left|h+\dd\psi\right|_g},\,\dd\eta \ket_g \,\dvg = 0.
\end{align*}
Integration by parts gives us
\begin{align*}
\int_M \eta\,\dd^*\left(\frac{ h+\dd\psi}{\,\left|h+\dd\psi\right|_g}\right)\,\dvg  = 0.
\end{align*}
As $\eta$ is arbitrary in the class $\{\eta \in C^1(M):\,\|\eta\|_{C^1} \leq c_0\}$, we obtain \eqref{EL}. \end{proof}

\smallskip
\noindent
{\bf 3.} With \emph{Claim B} at hand, let us conclude \emph{Claim A} by establishing the existence of minimiser to the functional $\mathcal{I}$ over $\mathcal{A}_h$. 
  
  \begin{proof}[Proof of Claim A]

We shall follow the direct method of calculus of variations.

First of all, note that $h \in \mathcal{A}_h$ with $\mathcal{I}[h]<\infty$, so $i_0:=\inf_{\alpha \in {\mathcal{A}_h}} \mathcal{I}[\alpha]$ is finite. Let $\{\alpha_j\}_{j \in \mathbb{N}} \subset {\mathcal{A}_h}$ be such that $\mathcal{I}[\alpha_j]= \|\alpha_j\|_{L^1(M,g)} \to i_0$ as $j \to \infty$. If $i_0=0$, then $\alpha_j \to 0$ (the zero 1-form) in the $L^1$-norm. But by assumption $h$ is not cohomologous to $0$. Thus we have $0<i_0<\infty$.

Now, by convergence of $\mathcal{I}[\alpha_j]$ we have that $\{\alpha_j\}$ is uniformly bounded in  $L^1$. In addition, it follows from the definition of ${\mathcal{A}_h}$ that $\|\alpha_j\|_{L^\infty(M,g)} \leq 10 \|h\|_{L^\infty(M,g)}$, so $\{\alpha_j\}$ is equi-integrable in $L^1$, so by the Dunford--Pettis theorem there exists a subsequence $\{\alpha_{j_k}\}$ $L^1$-weakly convergent to $\alpha_\star$. By lower semi-continuity of the $L^1$-norm with respect to the weak topology, we get
\begin{align*}
i_0 \leq \mathcal{I}\left[\alpha_\star\right] = \left\|\alpha_\star\right\|_{L^1(M,g)} \leq \liminf_{k \to \infty} \left\|\alpha_{j_k}\right\|_{L^1(M,g)} = i_0.
\end{align*}
In addition, the weak $L^1$-convergence implies that there is a further subsequence $\left\{\alpha_{j_{k_{\ell}}}\right\}$ converging \emph{a.e.} to $\alpha_\star$. So $\|\alpha_\star\|_{L^\infty(M,g)} \leq 10\|h\|_{L^\infty(M,g)}$. Also, the weak convergence $\alpha_{j_k}\weak \alpha_\star$ in $L^1$ implies convergence as 1-currents (\emph{i.e.}, in the sense of distributions), so $\alpha_\star$ is cohomologous to $h$  in the little regularity setting. Therefore, $\alpha_\star \in \mathcal{A}_h$. We can now conclude that $i_0:=\inf_{\alpha \in {\mathcal{A}_h}} \mathcal{I}[\alpha] =\mathcal{I}\left[\alpha_\star\right]$. \end{proof}

Finally, if the harmonic 1-form $h$ is everywhere non-vanishing, we may repeat the above arguments with the admissible class $\mathcal{A}_h$ (defined in \eqref{Ah}) replaced by 
\begin{equation}
\mathcal{A}_h' := \left\{ \beta \in [[h]]:\, \frac{1}{10} {\rm ess\,inf}_{M} |h|_g \leq |\beta|_g \leq 10 \|h\|_{L^\infty(M,g)}\text{ a.e.} \right\}.
\end{equation}
All the arguments go through; in particular, by the \emph{a.e.} convergence of a  subsequence, the limiting form $\alpha_\star$ is in $\mathcal{A}_h'$, hence  $\alpha_\star \neq 0$ \emph{a.e.}.

Thus, one may choose a nowhere vanishing representative of $\alpha_\star$ to conclude the proof.  \end{proof}

\end{appendices}

\bigskip
\noindent
{\bf Acknowledgement}.
This work has been done during SL's stay as a CRM--ISM postdoctoral fellow at  the Centre de Recherches Math\'{e}matiques, Universit\'{e} de Montr\'{e}al, McGill University, and Concordia University. SL would like to thank these institutions for their hospitality. SL also thanks the Shanghai Frontier Research Institute for Modern Analysis for its support during the finalisation of the paper.

MS is indebted to Vincent Borrelli and Amit Acharya for many helpful discussions. MS was supported in part by Simons Collaborative Research Grant 232531.

\end{document}